\documentclass[12pt,a4paper]{amsart}
\usepackage{amsfonts, amssymb, amsmath, amsthm}
\usepackage{graphicx} 
\usepackage{amscd}
\usepackage{color}
\usepackage{float}
\usepackage{paralist}
\usepackage[margin=5mm]{caption}
\usepackage{setspace}
\newtheorem{theorem}{Theorem}
\theoremstyle{plain}

\newtheorem{definition}[theorem]{Definition}
\newtheorem{corollary}[theorem]{Corollary}
\newtheorem{lemma}[theorem]{Lemma}\newtheorem{remark}{Remark}

\newtheorem{proposition}[theorem]{Proposition}
\newtheorem{mydef}[theorem]{Definition}
\def\bdefi{\begin{mydef}}\def\edefi{\end{mydef}}
\theoremstyle{definition}\newtheorem{example}{Example}

\numberwithin{equation}{section}

\newcommand{\secref}[1]{\S\ref{#1}}

\def\eps{\varepsilon}

\def\ds{\displaystyle}
\def\ri{{\rm i}}\def\Phibf{{\mathbf \Phi}}
\def\Psibf{{\mathbf \Psi}}\def\Gbf{{\mathbf G}}

\newcommand{\R}{\mathbb{R}}
\newcommand{\C}{\mathbb{C}}
\newcommand{\Q}{\mathbb{Q}}
\newcommand{\B}{\mathbb{B}}
\renewcommand{\P}{{\mathbb P}}

\newcommand{\N}{\mathbb{N}}

\newcommand{\Z}{\mathbb{Z}}

\newcommand{\cA}{{\mathcal A}}
\newcommand{\cAt}{\tilde{\mathcal A}}

\newcommand{\cF}{{\mathcal F}}
\newcommand{\cG}{{\mathcal G}}

\newcommand{\cL}{{\mathcal L}}

\newcommand{\cN}{{\mathcal N}}

\newcommand{\cS}{{\mathcal S}}

\newcommand{\cX}{{\mathcal X}}

\renewcommand{\phi}{\varphi}

\newcommand{\pa}{\partial}

\newcommand{\bA}{{\mathbf A}}
\newcommand{\bB}{{\mathbf B}}
\newcommand{\bb}{{\mathbf b}}
\newcommand{\bv}{{\mathbf v}}
\newcommand{\be}{{\mathbf e}}
\newcommand{\bR}{{\mathbf R}}
\newcommand{\bF}{{\mathbf F}}

\newcommand{\bW}{{\mathbf W}}
\newcommand{\bJ}{{\mathbf J}}

\newcommand{\phising}{\phi_{\text{sing}}}
\newcommand{\Vrev}{V_\text{rev}}

\DeclareMathOperator*{\essinf}{ess\,inf}

\DeclareMathOperator{\sign}{sign}

\newcommand{\bspm}{\left(\begin{smallmatrix}}\newcommand{\espm}{\end{smallmatrix}\right)}
\newcommand{\bpm}{\begin{pmatrix}}\newcommand{\epm}{\end{pmatrix}}

\def\blem{\begin{lemma}}\def\elem{\end{lemma}}
\def\bthm{\begin{theorem}}\def\ethm{\end{theorem}}
\def\bprop{\begin{proposition}}\def\eprop{\end{proposition}}
\def\bcor{\begin{corollary}}\def\ecor{\end{corollary}}
\def\beq{\begin{equation}}\def\eeq{\end{equation}}
\def\bpf{\begin{proof}}\def\epf{\end{proof}}
\def\bex{\begin{example}}\def\eex{\end{example}}

\def\bi{\begin{itemize}}\def\ei{\end{itemize}}
\def\ig{\includegraphics}\def\om{\omega}\def\sig{\sigma}
\def\Om{\Omega}\def\er{{\rm e}}\def\ri{{\rm i}}
\def\huga#1{\begin{gather} #1 \end{gather}}

\newcommand{\reff}[1]{(\ref{#1})}
\def\brem{\begin{remark}}\def\erem{\end{remark}}

\newcommand{\bcen}{\begin{compactenum}}\newcommand{\ecen}{\end{compactenum}}
\begin{document}
\title[Nonlinear Bloch Waves in the Gross-Pitaevskii Equation]{Bifurcation of Nonlinear Bloch Waves from the Spectrum in the Gross-Pitaevskii Equation} 

\author{Tom\'{a}\v{s} Dohnal}
\address{T. Dohnal \hfill\break 
Department of Mathematics, Technical University Dortmund\hfill\break
D-44221 Dortmund, Germany}
\email{tomas.dohnal@math.tu-dortmund.de}

\author{Hannes Uecker}
\address{H. Uecker \hfill\break 
Institute of Mathematics, Carl von Ossietzky University Oldenburg\hfill\break
D-26111 Oldenburg, Germany}
\email{hannes.uecker@uni-oldenburg.de}

\date{\today}


\begin{abstract} We rigorously analyze the bifurcation of stationary so called 
nonlinear Bloch waves (NLBs) from the spectrum in the Gross-Pitaevskii 
(GP) equation with a periodic potential, in arbitrary space dimensions. These 
are solutions which can be expressed as finite sums of 
quasi-periodic functions, and which in a formal asymptotic expansion 
are obtained from solutions of the so called 
algebraic coupled mode equations. Here we justify this expansion 
by proving the existence of NLBs and estimating the error 
of the formal asymptotics. The analysis is illustrated by numerical 
bifurcation diagrams, mostly in 2D. In addition, we illustrate some relations of NLBs to other classes of solutions of the GP equation, in 
particular to so called out--of--gap solitons and truncated NLBs, and 
present some numerical experiments concerning the stability of these 
solutions.
\end{abstract}

\maketitle

\section{Introduction}
 The nonlinear Schr\"odinger/Gross--Pitaevskii (GP) equation 
in $d\in\N$ dimensions, 
\huga{\label{tGP}
\ri\pa_t\psi=\Delta \psi
-V(x) \psi - \sigma |\psi|^2\psi, \quad x\in \R^d,\ t\in\R, 
}
with a real potential $V:\R^d\to \R$ 
is a canonical model in mathematics and physics. It appears in 
various contexts, e.g.,  
nonlinear optics \cite{SK02, Efrem_etal_2003}, and 
 Bose--Einstein conden\-sation \cite{losk03,aok06}. 
See also, e.g., \cite{sulem99, Y2010, fibich2015} for mathematical 
and modeling background. 
Plugging $\er^{\ri\om t}\phi(x)$ into \reff{tGP}, 
where  $\om/(2\pi)$ is the frequency of time--harmonic waves in 
nonlinear optics, 
and where  $\om$ is called the chemical potential in Bose--Einstein 
conden\-sation, 
 we obtain the stationary 
problem 
\beq\label{E:SNLS_dD} \omega \phi + \Delta \phi
-V(x) \phi - \sigma |\phi|^2\phi=0. 
\eeq 

Here we consider the case that 
the potential $V$ is real and periodic. For simplicity, we let $V$ be $2\pi-$periodic in each coordinate
direction, i.e., 
\[V(x+2\pi e_j) = V(x) \quad \text{ for all } x\in \R^d, j\in
\{1,\dots,d\},\] 
where $e_j$ denotes the $j$-th Euclidean unit vector
in $\R^d$. In other words, we consider the periodic lattice $2\pi\Z^d$. We make the basic assumption that 
$V \in H^{s-2}(\P)$ for some $s>\frac{d}{2}$, where $\P=(-\pi,\pi]^d$.
This smoothness assumption on $V$ ensures $H^s(\P)$-smoothness of 
linear Bloch waves, i.e., solutions of \eqref{E:SNLS_dD} with $\sigma=0$. 
See \S\ref{S:Bloch} for a review of spectral properties of 
$$
L=-\Delta +V
$$ 
and linear Bloch waves. 
For suitable $V$ the spectrum of $L$ shows so called spectral 
gaps and in recent years a focus has been on the bifurcation 
of so called gap solitons from the 
the zero solution at band edges into the gaps. These are localized solutions, which, in the near edge asymptotics have small amplitude and long 
wave modulated shape. In detail, the asymptotic expansion at $\om=\om_*+\eps^2\Om$ with $\Om=\pm 1$ is
\beq \label{as1} 
\phi(x)\sim\eps \sum_{j=1}^N A_j(\eps x)\xi_{n_j}(k^{(j)};x), 
\end{equation}
where $\xi_{n_j}(k^{(j)};\cdot)$, $j=1,\ldots,N$ are Bloch waves 
at the edge $\om_*$, 
and the $A_j$ are localized solutions of a system of  
(spatially homogeneous) nonlinear Schr\"odinger equations. 
See, for instance, \cite{sy07, DPS09, DU09,IW10}, and the references therein.  

Here we seek solutions of \eqref{E:SNLS_dD} which can be expressed as
finite a sum of $M$ quasi-periodic functions and call such solutions
\textit{nonlinear Bloch waves} (NLBs), with quasi--periodicities
determined from a selected finite subset of the Bloch waves at
$\omega$. NLBs have been studied in, for instance, \cite{DPS09,
  wyak09,zlw09,zw09,CP2012}, where in \cite{wyak09,zlw09,zw09} the
approaches are numerical and formal. They have been observed even
experimentally, see e.g. \cite{CMMCA02} for experiments in
Bose-Einstein condensates. In \cite{DPS09} the special case of a
bifurcation of NLBs into an asymptotically small spectral gap for a
separable periodic potential in two dimensions is studied
rigorously. In \cite{CP2012} the bifurcation of single component
($M=1$) NLBs in one dimension is proved, including results on 
secondary bifurcations and exchange of stability. 
Similarly to Bloch waves in
linear lattices NLBs can be understood as the fundamental bounded
oscillatory states of the nonlinear system. From the applied point of
view one motivation for studying NLBs is the continuation of gap--solitons to
``out-of-gap'' solitons, i.e., the continuation of localized solutions
from one band edge across the gap and into the spectrum on the other
side of the gap, where their tails start interacting with the NLBs.
For this reason, the study of bifurcation of NLBs from the
zero--solution has been mostly restricted to band edges.  Here we show
that nonlinear Bloch waves bifurcate in $\omega$ from {\em generic
  points in the spectrum} of $L$, and give their asymptotic expansions
in terms of solutions of the so called algebraic coupled mode
equations (ACME), together with error estimates.

\begin{figure}[ht]
  \begin{center}
\begin{tabular}{cc}
(a)&(b)\\
\ig[width=45mm]{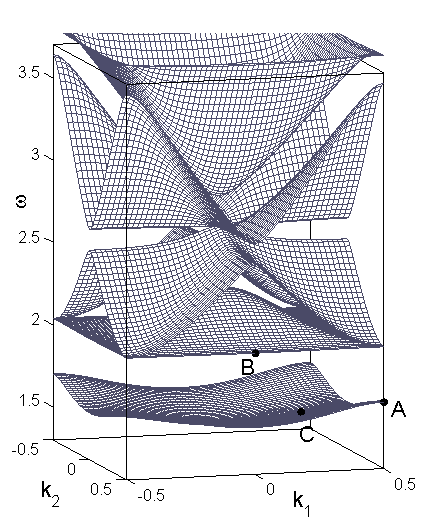}&
\ig[width=45mm, height=55mm]{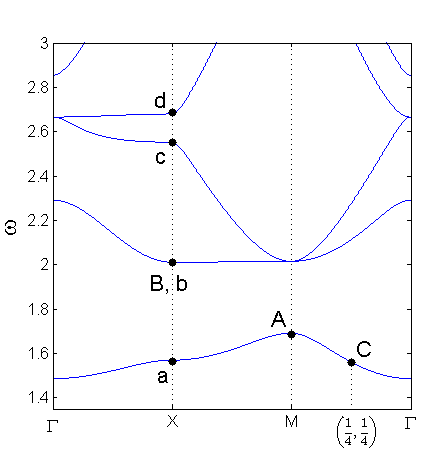}
\end{tabular}
  \end{center}
  \caption{{\small (a): Band structure of $L$ over the Brillouin zone $\B$
 for the periodic potential \eqref{E:V}; (b):
along the boundary $\Gamma-X-M-\Gamma$ of
    the, so called, irreducible Brillouin zone.}
}
  \label{F:band_str}
\end{figure}

In addition to the rigorous analysis we illustrate our results 
 numerically. For this we focus on 2D, as this is much richer than 
1D, and use the same potential as in \cite{DU09}, i.e.
\begin{equation}\label{E:V}
  V(x)=1+4.35 W(x_1)W(x_2), \quad x \in [-\pi,\pi]^2
\end{equation}
with
\[W(s)=\frac{1}{2}\left[\tanh\left(7\left(s+\frac{3\pi}{5}\right)\right)+\tanh\left(7\left(\frac{3\pi}{5}-s\right)\right)\right].\]
This represents a square geometry with smoothed-out edges. The function in \eqref{E:V} is extended periodically to $\R^2$ to obtain $V:\R^2\to \R$.
The numerical band structure of $L$ over the 
Brillouin zone  $\B := (-1/2,1/2]^d$, 
and also along the boundary of the irreducible Brillouin zone, 
is plotted in Fig.~\ref{F:band_str}(a),(b), respectively. We denote the so called high symmetry points in $\B$ for $d=2$ by
$$\Gamma:=(0,0), X:=(1/2,0), X':=(0,1/2), \quad \text{and } M:=(1/2,1/2).$$

\begin{example}\label{Ex:NLB_intro}
Figure \ref{f2} shows a 
numerical bifurcation diagram of single component ($M=1$) NLBs for $k=X$, calculated 
with the package \texttt{pde2path}
\cite{p2p,p2p2}, together with example plots on the bifurcating 
branches. A branch of NLBs bifurcates from 
the zero solution at $\omega=\omega_*$ for any $\omega_*$ attained by one of the band functions at $k=X$, i.e. at the $\omega-$coordinate of any of the points $a,b,c,d$ in Fig. \ref{F:band_str} (b). See  \S\ref{S:Bloch} for the definition of band functions.
\begin{figure}[ht]
  \begin{center}
\ig[width=60mm,height=60mm]{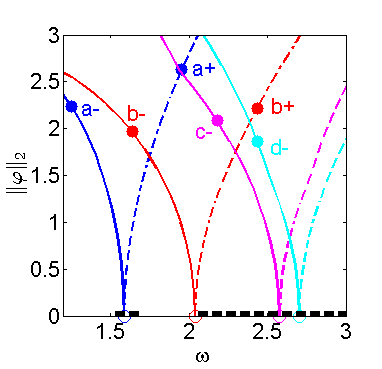}\hspace{10mm}
\raisebox{35mm}{\begin{minipage}{75mm}
\ig[width=35mm]{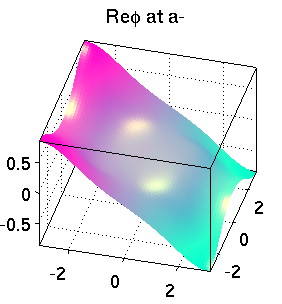}\ig[width=35mm]{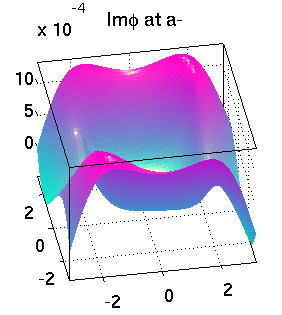}
\ig[width=35mm]{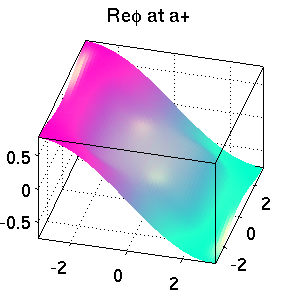}\ig[width=35mm]{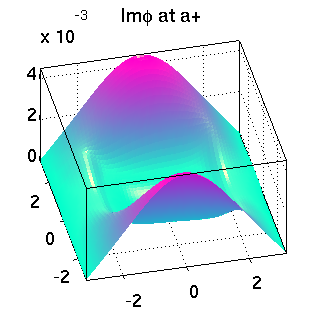}
\end{minipage}}

\ig[width=35mm]{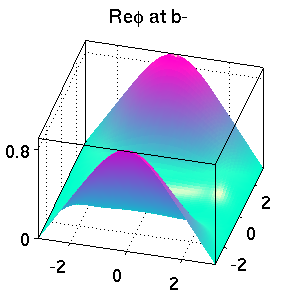}\ig[width=35mm]{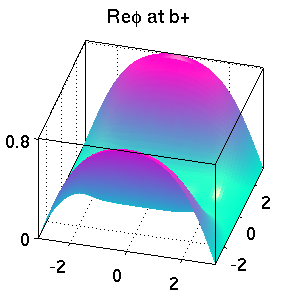}
\ig[width=35mm]{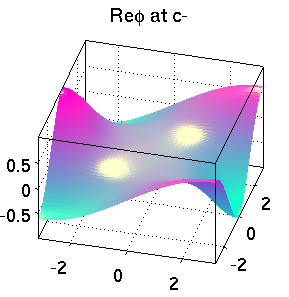}\ig[width=35mm]{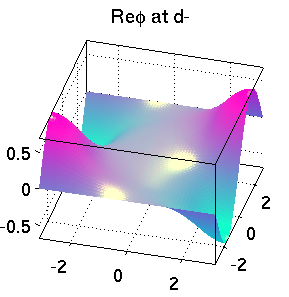}
  \end{center}
 \caption{{\small Example \ref{Ex:NLB_intro}: Bifurcation diagram of the first four bifurcating branches for $k=X$, i.e. branches bifurcating from points a-d in Fig. \ref{F:band_str} (b). Spectral bands are indicated by the black dashed line. The sign $\pm$ in the branch labels stands 
for $\sig=\pm 1$. Small panels: 
example solution plots of NLBs from the bifurcation 
diagram, over the fundamental cell $x\in (-\pi,\pi)^2$. 
At bifurcation we choose a real Bloch wave. Then the imaginary 
parts are small near bifurcation, and we only plot them for 
a$\pm$. Roughly horizontal axis corresponds to $x_1$ in all plots. }
}\label{f2}
\end{figure}

In \secref{S:numerics} we explain the 
method behind Fig.~\ref{f2}, and study in detail the bifurcations of NLBs 
at the points marked A,B,C in Fig.~\ref{F:band_str}(b), 
relating the numerical calculations to our analysis. 
\end{example}

As already said, one motivation for studying NLBs are the intriguing properties
of their interaction with localized solutions, which we illustrate 
numerically in \S\ref{ogssec}. 
 For instance, when a
gap soliton is continued from the gap into the spectrum, we get a so
called ``out--of--gap'' soliton (OGS) with oscillating (delocalized)
tails, see also \cite{ys03,JKKK11}. In 1D, numerically these OGS can
be seen to be homoclinic orbits approaching NLBs, 
and essentially the same happens in 2D. 
Moreover, the NLB can form building blocks of so called truncated NLBs
(tNLBs), see also \cite{aok06,wyak09}.  These are localized solutions
for $\om$ in a gap which are close to a NLB on some finite interval
but approach $0$ as $|x|\to\infty$. Then, continuing a branch of tNLBs 
from the gap into the spectrum, the same
interaction scenario as for GS happens, i.e., the tails of 
the tNLBs pick up NLBs bifurcating from the gap edge into 
the spectrum, and the tNLBs become delocalized, for which 
we use the acronym dtNLB. Note that both gap solitons and
tNLBs have been observed experimentally, see
e.g. \cite{BOP2012} for experiments in optical lattices. However, even
in 1D at present it is unclear how to analyze OGS, tNLBs and dtNLBs 
rigorously, i.e., so far there
only exist heuristic asymptotics, see \S\ref{ogssec} for further
comments. 

Stability of most of these solutions is an open problem. 
Thus, at the end of the paper 
we also give a numerical outlook on this, and obtain stability 
of some NLBs in 1D, and, consequently, stability of some tNLBs and 
some OGS and dtNLBs. In 2D, we did not find stable NLBs for the potential 
$V$ from \reff{E:V}, and we close with summarizing the open questions. 
The broad spectrum of applications of NLBs clearly motivates our rigorous bifurcation analysis.

In the remainder of this introduction we explain the linear band
structure, a simple analytical bifurcation result, formulate the main
theorem, and describe the structure of the paper in more detail.
 
\subsection{Linear Bloch waves}\label{S:Bloch}
For $k$ in the Brillouin zone, $k\in\B := (-1/2,1/2]^d$, consider the Bloch eigenvalue problem
\beq\label{E:Bloch_ev_prob}
\begin{aligned}
(-\Delta + V(x)) \xi_n(x,k) &= \omega_n(k) \xi_n(x,k),\quad x\in \P\\
\xi_n(x+2\pi e_m,k) &= e^{2\pi\ri k_m}\xi_n(x,k), \quad m \in \{1,\dots,d\}, 
\end{aligned}
\eeq 
where $e_m$ is the $m-$th Euclidean unit vector in $\R^d$. 
The spectrum of $L=-\Delta +V$ is continuous and is given by the
union of the bands defined by the band structure $(\omega_n(k))_{n\in
  \N}$, i.e.
\[\sigma(L)=\bigcup_{\substack{n\in \N \\ k\in
    \B}}\omega_n(k) = \bigcup_{l\in \N}[s_{2l-1},s_{2l}], \ \text{
  where } s_{2l-1} < s_{2l} \leq s_{2l+1} \ \text{ for all } l\in
\N,\] 
see Theorem 6.5.1 in \cite{Eastham}. The functions $k\mapsto
\omega_n(k)$ are called band functions. The Bloch waves $\xi_n(x,k)$
have the form $\xi_n(x,k)=p_n(x,k)e^{\ri k\cdot x}$ with $p_n(x+2\pi
e_m,k) = p_n(x,k)$ for all $m \in \{1,\dots,d\}$ and all $x\in
\R^d$. Clearly, both $\omega_n(k)$ and $\xi_n(x,k)$ are $1-$periodic in each component of $k$.  We assume the normalization
$$\|\xi_n(\cdot, k)\|_{L^2(\P)}=\|p_n(\cdot, k)\|_{L^2(\P)}=1 \quad \forall n \in \N \ \forall k \in \B.$$ 

For a given point $(k,\omega)\in\B\times \R$ in the band
structure, i.e. with $\omega =\omega_n(k)$ for some $n\in \N$, also
the point $(-k,\omega)$ lies in the band structure, which follows from
the symmetry 
\beq\label{E_sym_om} 
\omega_n(k)=\omega_n(-k) \ \text{for all} \ n\in \N, k\in \B.  
\eeq 
	This symmetry is due to the
equivalence of complex conjugation and replacing $k\mapsto -k$ in the
eigenvalue problem \eqref{E:Bloch_ev_prob}. Hence, we also have the conjugation symmetry of the Bloch waves, namely
\beq\label{E:sym_Bloch} 
\xi_n(x,k) = \overline{\xi_n(x,-k)}.  
\eeq
For $k\in \partial \B\cap \B$ we have $-k\in \partial \B
\setminus \B$ and the point $-k$ must be understood as the
$\Z^d$-periodic image within $\B$. When $k$ is one of the so called
\textit{high symmetry points}, i.e. $k_m \in \{0,1/2\}$ for all $m\in \{1,
\dots,d\}$, then $k$ and $-k$ are identified via this
periodicity. Equation \eqref{E:sym_Bloch} then implies that
$\xi_n(x,k)$ is real. This can be seen directly from the eigenvalue
problem \eqref{E:Bloch_ev_prob}, where $k_m \in \{0,1/2\}$ for all
$m\in \{1, \dots,d\}$ implies that the boundary condition is real such
that a real eigenfunction must exist. Note that the above $k-$symmetries rely only on the realness of $V$.

\subsection{The Bifurcation Problem}

\begin{remark}\label{Rem:real}{
In the simplest scenario we can look for real solutions of
\eqref{E:SNLS_dD} with the quasi-periodic boundary conditions given by
a single vector $k_*\in \B$, i.e.  
\beq\label{E:qp_kstar} 
\phi(x+2\pi e_m) = \phi(x) e^{2\pi \ri k_{*,m}} \quad \text{for all } x\in \R^d,
m\in \{1,\dots, d\}.  
\eeq 
In this case the realness condition on $\phi$ requires
\beq\label{E:kstar_cond}
k_*\in \left\{0,\tfrac{1}{2}\right\}^d,
\eeq
such that the sought solution is $2\pi$-periodic or $2\pi$-antiperiodic in each coordinate direction. We study bifurcations in the parameter
$\omega$. Classical theory for bifurcations at simple eigenvalues, e.g., Theorem 3.2.2 in \cite{Nirenb_1974}, shows that if $\omega_* =
\omega_{n_*}(k_*)$ for exactly one $n_*\in\N$, i.e. $\omega_*$ is a
simple eigenvalue of $L$ under the boundary conditions \eqref{E:qp_kstar}, then $\omega=\omega_*$ is a bifurcation
point. To this end define $f(\phi,\omega)=\omega \phi +\Delta \phi -V(x)\phi-\sigma \phi^3$ and study $f(\phi,\omega)=0$ on $\P$ under the boundary conditions  \eqref{E:qp_kstar}.
We have $f(0,\omega)=0$ for all $\omega\in \R$ and $f_\phi(0,\omega)=\omega-L$. As $\omega_*$ is a simple eigenvalue, we have the one dimensional kernel
$$\mbox{Ker}(f_\phi(0,\omega_*))= \xi_{n_*}(x,k_*)\R.$$
Because $L$ with \eqref{E:qp_kstar} and \eqref{E:kstar_cond} is self adjoint, we have $\mbox{ Ran}(f_\phi(0,\omega_*))\perp_{L^2(\P)}\text{Ker}(f_\phi(0,\omega_*))$.
The transversality condition $f_{\omega\phi}(0,\omega_*)\xi_{n_*}(x,k_*) \notin \text{Ran}(f_\phi(0,\omega_*))$ of Theorem 3.2.2 in \cite{Nirenb_1974} thus holds because $f_{\omega\phi}(0,\omega_*)\xi_{n_*}(x,k_*)=\xi_{n_*}(x,k_*)\perp_{L^2(\P)}\text{Ran}(f_\phi(0,\omega_*))$. As a result, the theorem guarantees the existence of a unique non-trivial branch of solutions bifurcating from $\omega=\omega_*$.
}
\end{remark}

\begin{remark} {
Without the restriction to real solutions the eigenvalue $\omega_*$ is never simple due to invariances. In the real variables $\Phi:=(\phi_R,\phi_I)^T$, where $\phi=\phi_{R}+\ri \phi_I$, the problem becomes
$$\cG(\Phi,\omega)=\begin{pmatrix}\omega \phi_R + \Delta\phi_R - V(x) \phi_R-\sigma(\phi_R^2+\phi_I^2)\phi_R \\
\omega\phi_I + \Delta\phi_I - V(x)\phi_I -\sigma(\phi_R^2+\phi_I^2)\phi_I
\end{pmatrix}=0.$$
Since \eqref{E:SNLS_dD} possesses the phase invariance and the complex conjugation invariance, we get that $\cG$ is $O(2)$ invariant, i.e. 
$$\cG(\gamma \Phi,\omega)=\gamma\cG(\Phi,\omega) \text{ for all } \gamma \in \Gamma:=\left\{\bspm 1 & 0 \\ 0 & -1\espm,\bspm \cos \theta & -\sin \theta \\ \sin\theta & \cos\theta\espm: \theta\in [0,2\pi)\right\}.$$
Bifurcations can now be studied using the equivariant branching lemma, see  e.g. \cite[\S5]{mei2000}, by restricting to a fixed point subspace of a subgroup of $\Gamma$. The only nontrivial subgroup is $\{\bspm 1 & 0 \\ 0 &1\espm,\bspm 1 & 0 \\ 0 & -1\espm\}$ with the fixed point subspace being the vectors $\Phi$ with $\phi_I=0$ corresponding to real solutions of \eqref{E:SNLS_dD}. Therefore, this leads again to real solutions. Nevertheless, more complicated solutions than the single component ones in Remark \ref{Rem:real} can be studied. The most general real ansatz is
\beq\label{E:ans_ebl}
\phi(x) = \sum_{j=1}^{2q+r}\phi_j(x), \quad \phi_j(x+2\pi e_m) = e^{\ri 2\pi k_{m}^{(j)}}\phi_j(x), m=1,\dots,d
\eeq
with $q,r\in \N_0$, with $k^{(j)}\in \B$ for all $j=1,\dots,2q+r$, such that $k^{(j+q)}\dot{=}-k^{(j)}, \phi_{j+q}=\overline{\phi_j}$ for all $j=1,\dots,q$, and with $k^{(j)}\in \{0,1/2\}^d, \phi_j(x)\in \R$ for $j=2q+1,\dots,2q+r$. Here $\dot{=}$ means equality modulo 1 in each coordinate. While the use of the equivariant branching lemma should describe the bifurcation problem and produce the effective Lyapunov-Schmidt reduction, we choose to carry out a detailed analysis without this tool in order to obtain more explicit results. This will allow us to provide estimates of the asymptotic approximation error. 
}
\end{remark}

Our aim is to prove a general bifurcation theorem for NLBs, and, 
moreover, to derive and justify an effective asymptotic model related to the Lyapunov-Schmidt reduction of the bifurcation problem including an estimate on the asymptotic error. In our approach we select a frequency $\omega_*$ in the spectrum and choose $N$ points $\{k_*^{(1)},\dots, k_*^{(N)}\}\subset \B$ in the level set of the band structure at $\omega_*$, such that for each $j$ we have $\omega_*=\omega_{n_j}(k_*^{(j)})$ for some $n_j\in \N$.  Our method requires that each of the points $\{k_*^{(1)},\dots, k_*^{(N)}\}$ is either one of the high symmetry points $k\in \{0,1/2\}^d$ or belongs to a pair $k,l$ with $l\dot{=}-k$.  See (H1)--(H6) on page \pageref{kdef} for a summary 
of our assumptions. 
We seek NLBs bifurcating from $\omega_*$ and having the asymptotic form 
\beq\label{E:asymp-intro} 
\phi(x) \sim \eps \sum_{j=1}^N A_j \xi_{n_j}(x,k^{(j)}), 
\eeq at $\omega=\omega_*\pm\eps^2$. The coefficients $A_j$, i.e. the
(complex) amplitudes of the waves, are given by solving the ACME as an effective algebraic system of $N$ equations. 
 Generally a sum of $N$ quasiperiodic functions with the quasiperiodicity of each given by one of the vectors $k^{(j)}$ cannot be an exact solution of \eqref{E:SNLS_dD} as the nonlinearity generates functions with other quasiperiodicities. Our ansatz for the exact solution is thus
$$\phi(x) =\sum_{j=1}^M \phi_j(x), \quad \phi_j(x+2\pi e_m)= e^{\ri 2\pi k^{(j)}_m}\phi_j(x), \quad m=1, \dots,d$$
with $M\geq N$ and $k^{(j)}=k_*^{(j)}$ for $j=1,\dots,N$. Importantly, this set $\{k^{(1)},\dots, k^{(M)}\}$ (defined in \eqref{kdef}) can be a proper
subset of the level set. The subset may be finite even if the level
set is, for instance, uncountable. In fact, our assumption (H4) ensures the finiteness.
Besides, the subset $\{k^{(1)}, \dots, k^{(N)}\}$ can be much
smaller than $\{k^{(1)}, \dots, k^{(M)}\}$ and hence the effective
ACME-system can be rather small. The subset has to satisfy only (H2-H6).

The major assumptions of our analysis are rationality (assumption (H4))
and certain non-resonance conditions (H5) on the $k$-vectors
$\{k^{(1)}, \dots, k^{(N)}\}$. In addition, the solutions of the
coupled mode equations need to satisfy certain symmetry
(``reversibility'') and non-degeneracy conditions, see
Definitions \ref{D:nondegen} and \ref{D:reverse}, in order for us to guarantee that
\eqref{E:asymp-intro} approximates a solution $\phi$ of
\eqref{E:SNLS_dD}.The main result is the following 
\bthm\label{T:main} 
Assume (H1)-(H6). There exist $\eps_0>0$ and $C>0$ such that 
for all $\eps\in(0,\eps_0)$ the following holds. If 
the  ACMEs \eqref{E:CME_NLB} have a reversible non-degenerate 
solution $\bA\in
\Vrev$, then \eqref{E:SNLS_dD} with $\omega =\omega_* +\eps^2 \Omega$
has a nonlinear Bloch wave solution $\phi$ of the form
\eqref{E:NLB_ans}, and 
\[\left\|\phi(\cdot)-\eps\sum_{j=1}^N
  A_j\xi_{n_j}(\cdot,k^{(j)})\right\|_{H^s(\P)}\leq C \eps^3.\]
 \ethm

 There are three relatively straightforward generalizations of the
 result. Firstly, the periodic lattice $2\pi\Z^d$ can be replaced by
 any lattice $\{\sum_{j=1}^dm_ja_j:m\in \Z^d\}$ with linearly
 independent vectors $\{a_1,\dots,a_d\}\subset \R^d$. Of course, the
 periodicity cell $\P$ and the Brillouin zone $\B$ have to be
 redefined accordingly. Except for the examples in
 \S\ref{S:ACMEs_N12} the results, in particular Theorem
 \ref{T:main}, hold for a general lattice. Secondly, the nonlinearity
 $|\phi|^2\phi$ can be replaced by other locally Lipschitz
 nonlinearities $f(\phi)$ which are phase invariant and satisfy
 $f(\phi)=o(\phi)$ for $\phi \to 0$. This may, however, change the
 powers of $\eps$ in the expansion and the error estimate. Also, the
 linear operator $L$ can be generalized to self adjoint second order
 differential operators with periodic coefficients such that the
 asymptotic distribution of eigenvalues $\omega_n(k)$ remains that in
 \eqref{E:as_bands}.

\subsection{The Structure of the Paper}\label{S:struct}
In \S\ref{S:formal} we present 
a formal asymptotic approximation of nonlinear Bloch
waves and a derivation of the ACMEs as effective amplitude equations. 
In \S\ref{S:ass_LS} we
pose conditions on the solution ansatz and the band structure which
are necessary for our analysis, and apply the Lyapunov-Schmidt
decomposition to the bifurcation problem. The invertible part of the decomposition is estimated in 
\S\ref{S:reg}. The singular part and its relation to the ACMEs is
described in \S\ref{S:sing}, where also the proof of the main theorem is completed. In \S\ref{S:ACMEs_N12} we present the
ACMEs and their solutions in the scalar case ($N=1$) and in the case
of two equations ($N=2$).  
Section \ref{S:numerics} presents
numerical computations of nonlinear Bloch waves in two dimensions
$d=2$ for $N=1$ and $N=2$. The convergence rate of the approximation
error is confirmed by numerical tests. 
Finally, in \S\ref{ogssec} we give a numerical outlook on the 
interaction of localized solutions with NLBs, first for some 1D and 2D 
GS, and second for tNLBs, and we report numerical experiments on stability of 
NLBs and other solutions. 


\section{Formal Asymptotics}\label{S:formal}
Let $\omega_*\in \sigma(L)$ and choose $N\in \N$
vectors $k^{(1)},\dots, k^{(N)}\in \B$ in the level set of the band
structure at $\omega_*$.  For the asymptotics of nonlinear Bloch waves
near $\omega_*$ we make an analogous ansatz to that used in
\cite[\S3]{DPS09} for nonlinear Bloch waves near band edges in \reff{E:SNLS_dD} 
with a {\em separable} periodic potential.  Formally we write
\beq\label{E:ansatz_mult} \phi(x)\sim \eps \sum_{j=1}^N A_j
\xi_{n_j}(x,k_*^{(j)}) +\eps^3\sum_{j=1}^N\phi^{(1)}_j(x) \quad
\text{for} \ \omega =\omega_*+\eps^2\Omega \quad (\eps \to 0), \eeq
where the amplitudes $A_j\in \C$ are to be determined and where
$\phi^{(1)}_j$ satisfies the quasiperiodicity given by the vector
$k_*^{(j)}$.

Substituting \eqref{E:ansatz_mult} in \eqref{E:SNLS_dD} 
we get at $O(\eps^3)$ for each $j\in \{1,\dots,N\}$
$$(-\Delta +V(x)-\omega_*)\phi^{(1)}_j(x) = \Omega A_j \xi_{n_j}(x,k_*^{(j)}) -\sigma \sum_{(\alpha,\beta,\gamma) \in \cA_j}A_\alpha \overline{A}_\beta A_\gamma \xi_{n_\alpha}(x,k_*^{(\alpha)})\overline{\xi_{n_\beta}(x,k_*^{(\beta)})} \xi_{n_\gamma}(x,k_*^{\gamma}),$$
where
\huga{\label{ajdef}
\cA_j=\{(\alpha,\beta,\gamma)\in \{1,\dots, N\}^3: \
k_*^{(\alpha)}-k_*^{(\beta)}+k_*^{(\gamma)}-k_*^{(j)}\in \Z^d\}.
}
The condition $(\alpha,\beta,\gamma)\in \cA_j$ in the sum
ensures that the nonlinear terms have the same quasi-periodicity as
$\phi^{(1)}_j$.  Nonlinear terms generated by the ansatz \eqref{E:ansatz_mult} and having other
quasi-periodicity than one of those defined by $k^{(j)}_*,
j=1,\dots,N$ have been ignored in this formal calculation.

Imposing the solvability condition, i.e.~making the right hand side
$L^2$-orthogonal to $\xi_{n_j}(\cdot,k_*^{(j)})$, we get the
\textit{algebraic coupled mode equations} (ACMEs)
\beq\label{E:CME_NLB} \Omega A_j - \cN_j(A_1,\dots, A_N)=0, \qquad
j\in \{1,\dots, N\}, \eeq
\begin{align}
\cN_j = & \sigma \sum_{(\alpha,\beta,\gamma)\in \cA_j} \mu_{\alpha,\beta,\gamma,j} A_\alpha\overline{A_\beta}A_\gamma, \notag \\
\mu_{\alpha,\beta,\gamma,j} = & \int_{\P}\xi_{n_\alpha}(x,k_*^{(\alpha)})\overline{\xi_{n_\beta}(x,k_*^{(\beta)})}\xi_{n_\gamma}(x,k_*^{(\gamma)})\overline{\xi_{n_j}(x,k_*^{(j)})}dx. \label{E:mu}
\end{align}

To make the approximation \eqref{E:ansatz_mult} rigorous, we must account for the nonlinear terms left out above and provide an estimate on the correction $\phi(x)-\eps \sum_{j=1}^N A_j \xi_{n_j}(x,k_*^{(j)})$.

\section{Solution Ansatz, Assumptions, Lyapunov-Schmidt Decomposition}\label{S:ass_LS}
As mentioned above, one of the difficulties of the analysis is that
for a sum of $N$ functions $f_1, \dots, f_N$ with distinct
quasi-periodic conditions the nonlinearity $|f_1+\dots
+f_N|^2(f_1+\dots +f_N)$ can generate functions with a new
quasi-periodicity. If the $k$-points defining these new
quasi-periodic boundary conditions lie in the $\omega_*$-level set of
the band structure, then a resonance with the kernel of the linear
operator occurs. Also, if the points generated by a repeated iteration
of the nonlinearity merely converge to the level set, our techniques
fail because a lower bound on the inverse of the linear operator
cannot be obtained. These obstacles are avoided if for a selected
$\omega_*\in \sigma(L)$ assumptions (H4) and (H5) below
hold. 

We select $N$ points $\{k_*^{(1)},\dots, k_*^{(N)}\}\subset \B$ in the $\omega_*$-level set of the band structure. Suppose we
seek solutions of \eqref{E:SNLS_dD} with $\phi$ given by the sum of quasiperiodic functions. The
ansatz $\phi(x) =\sum_{j=1}^N \phi_j(x)$ with quasiperiodic $\phi_j$
such that $\phi_j(x+2\pi e_m)=e^{\ri 2\pi k_{*,m}^{(j)}}\phi_j(x)$ for
all $x\in \R^d, m\in \{1,\dots,d\}$ can be a solution of
\eqref{E:SNLS_dD} only if each term generated by the nonlinearity
applied to this sum has quasiperiodicity defined by one of the vectors
in $\{k_*^{(1)},\dots,k_*^{(N)}\}$, i.e. if the \textit{consistency
  condition} 
\beq\label{E:consist} S_3(\{k_*^{(1)},\dots,
k_*^{(N)}\})\subset \{k_*^{(1)},\dots, k_*^{(N)}\} + \Z^d, 
\eeq
where
$$S_3: \{k_*^{(1)},\dots,k_*^{(N)}\} \to \{k_*^{(\alpha)}-k_*^{(\beta)}+k_*^{(\gamma)}: 1\leq \alpha,\beta,\gamma\leq N\},$$
is satisfied. In other words the consistency condition \eqref{E:consist}
says that all combinations $(\alpha,\beta,\gamma)$ for
$\alpha,\beta,\gamma\in \{1,\dots,N\}$ must lie in
$\cup_{j=1}^N\cA_j$, with $\cA_j$ from \reff{ajdef}. 

An example of a consistent ansatz for $N>1$ is $N=2, d=2$ with
$k_*^{(1)}=X=(1/2,0), k_*^{(2)}=X'=(0,1/2)$, like e.g. for
$\omega_*=s_3$ in \cite{DU09}. On the other hand, for $\omega_*=s_5$,
where $N=4, k_*^{(1)}=(k_c,k_c), k_*^{(2)}=(-k_c,k_c),
k_*^{(3)}=(-k_c,-k_c), k_*^{(4)}=(k_c,-k_c)$ with $k_c \approx 0.439$,
see Sec. 3.2.2.5 in \cite{DU09}, the ansatz is
inconsistent. It is also inconsistent for typical $\omega_*$ in the interior of $\sigma(L)$ with generic $\{k_*^{(1)},\dots,k_*^{(N)}\}$ in the level set. Therefore, we drop the consistency condition and pursue
the more general case where the nonlinearity generates quasiperiodic
functions with quasi-periodicity vectors $k$ not necessarily contained
in $\{k_*^{(1)},\dots, k_*^{(N)}\}$. 

Hence, we define the set of $k$-points generated by iterations of the nonlinear operator
\begin{align}
  & K := \{k\in \B : k\in S_3^p(\{k_*^{(1)},\dots,k_*^{(N)}\}) +
  \Z^d \text{ for some } p\in \N\}, \label{kdef}
\end{align}
and write, with $M\geq N$, 
\huga{\label{kdefb}
K=(k^{(j)})_{j=1}^M, \ \text{where} \ k^{(i)}=k_*^{(i)} \text{ for
} i=1,\dots,N.
}
At this point $M=\infty$ is possible but as explained below, our
assumption (H4) ensures $M<\infty$, i.e. only finitely many new vectors $k$ are
generated.  Thus we can search for a solution in the form of the sum of
finitely many quasiperiodic functions
\beq\label{E:NLB_ans} \phi(x) =\sum_{j=1}^M \phi_j(x), \qquad
\phi_j(x+2\pi e_m)=e^{\ri 2\pi k^{(j)}_m}\phi_j(x), \ m=1,\dots,d \eeq
with $\phi_j\in H^s(\P)$. The choice of the function space for $\phi_j$ is made clear below.

We make the following assumptions: 
\begin{itemize}
\item[{\bf (H1)}] 
$V \in H^{s-2}(\P)$ for some $s>\frac{d}{2}$, where $\P=(-\pi,\pi]^d$; 
\item[{\bf (H2)}] $\omega_*\in \sigma(L)$ and $k_*^{(1)}, \dots, k_*^{(N)}\in \B $ are
  points in the $\omega_*$-level set of the band structure, i.e., there are $n_1,\dots, n_N\in \N$ such that
$$\omega_{n_1}(k_*^{(1)})=\dots=\omega_{n_N}(k_*^{(N)})=\omega_*;$$ 
\item[{\bf (H3)}] each point $k_*^{(j)} \in \{k_*^{(1)}, ..., k_*^{(N)}\}$ is repeated according to the multiplicity of $\omega_*$ at $k=k_*^{(j)}$. In detail, if $q\geq 1$ band functions $\omega_{m_1}, ..., \omega_{m_q}$ touch at $(k,\omega)=(k_*^{(j)},\omega_*)$, then $q$ points in  $\{k_*^{(1)}, ..., k_*^{(N)}\}$ equal $k_*^{(j)}$ and $\{m_1, ..., m_q\} \subset \{n_1, ..., n_N\}$;
\item[{\bf (H4)}] the points $k_*^{(1)}, \dots, k_*^{(N)}\in \B $ have
  rational coordinates, i.e.
$$k_*^{(1)}, \dots, k_*^{(N)} \in \Q^d\cap \B;$$
\item[{\bf (H5)}] the intersection of the set $K$ with the level set of the band structure at $\omega=\omega_*$
  is exactly the set $\{k_*^{(1)}, \dots, k_*^{(N)}\}$, i.e.
$$K\cap \cL_{\omega_*}= \{k_*^{(1)}, \dots, k_*^{(N)}\},$$
where
$$\cL_{\omega_*} := \{k\in \B: \omega_n(k)=\omega_* \text{ for some }n\in \N\};$$
\item[{\bf (H6)}] for each $k_*^{(j)}\in \{k_*^{(1)}, \dots,
  k_*^{(N)}\}$ the reflection w.r.t. the origin lies in the set too, i.e.
$$k_*^{(j)}\in \{k_*^{(1)}, \dots, k_*^{(N)}\} \ \text{if and only if} \  k_*^{(j')}\in \{k_*^{(1)}, \dots, k_*^{(N)}\},$$
where $\B \ni k_*^{(j')}\dot{=}-k_*^{(j)}$ and $\dot{=}$ denotes
congruence with respect to the $1$-periodicity in each component.
\end{itemize}

With (H3) the bifurcation from multiple Bloch eigenvalues is
allowed. In one dimension ($d=1$) multiplicity is at most two, which
occurs for so called finite band potentials, see
e.g.~\cite{Maier_2008,D14}, and only at $k=0$ or $k=1/2$. In higher
dimensions ($d>1$) crossing or touching of band functions is abundant
in generic geometries (potentials $V$). 
Our results thus apply also to Dirac points in two dimensions studied,
e.g., in \cite{FW12}.

Due to the rationality condition (H4) the sought solution
\eqref{E:NLB_ans} is, in fact, periodic. (H4) also ensures that the
set $K$ is finite ($M<\infty$). Indeed, iterating the operator $S_3$
on a set of points with rational coordinates on a $d$-dimensional
torus generates a periodic orbit, i.e.~only finitely many distinct
points are generated, and the number $M$ depends solely on
$k_*^{(1)},\dots,k_*^{(N)}$. Condition (H4) is satisfied, e.g. if
$\{k_*^{(1)},\dots,k_*^{(N)}\}$ is a subset of the high symmetry
points of $\B$, i.e. $k_{*}^{(j)}\in \{0,1/2\}^d$ for all
$j=1,\dots,N$. This is frequently the case for the locations of
extrema defining a spectral edge. In general, (H4) is, however, a
serious limitation, and removing this assumption would be a major
improvement.

The non-resonance condition (H5) is satisfied, for instance, if
$\omega_* \in \partial \sigma(L)$, i.e. for $\omega_*$ at
one of the band edges, and $k_*^{(1)}, \dots, k_*^{(N)}$ are all the
extremal points of the band structure at which the edge $\omega_*$ is
attained.

The symmetry condition (H6) is needed in the persistence step of the
proof, see \S\ref{S:sing}. Note that if $k^{(j)}\in \cL_{\omega_*}$,
then also $k_*^{(j')}\in \cL_{\omega_*}$ by \eqref{E_sym_om} and the
periodicity in $k$. For $k^{(j)}\in \text{int}\B$, clearly,
$k_*^{(j')}=-k_*^{(j)}$. For $k^{(j)}\in \pa\B\cap \B$ is
$-k_*^{(j)}\in \pa\B\setminus \B$ (e.g. for $d=2, -k_*^{(j)}=(-1/2,a)$
with $a\in (-1/2,1/2)$) and then $k_*^{(j')}$ is the $\Z^d$-periodic
image of $-k_*^{(j)}$ within $\B$ (for the example
$k_*^{(j')}=(1/2,a)$). Moreover, also (H6) is automatically satisfied 
if $\{k_*^{(1)},\dots,k_*^{(N)}\}$ is a subset of the high symmetry
points of $\B$. 

Note again that $(k^{(i)})_{i=1}^N$ as well as $K$ may be proper subsets of $\cL_{\omega_*}$. This is, for instance, 
the case in 2D if $\cL_{\omega_*}=\{(1/2,0),(0,1/2)\}$,  where we could choose 
$N=1$ (if $\om_*$ is simple), 
and $K=\{(1/2,0)\}$ or $K=\{(0,1/2)\}$, which yields two decoupled 
scalar bifurcation problems; see \S\ref{S:N2} for further discussion.

The remaining two assumptions in Theorem \ref{T:main} are 
non-degeneracy and reversibility of $\bA$, defined as follows. 
\begin{definition}\label{D:nondegen}
  $\bA \in \C^N$ is a \textit{non-degenerate} solution of
  \eqref{E:CME_NLB} if the Jacobian\footnote{Strictly speaking, the
    problem should first be rewritten in real variables to define a
    Jacobian, see the discussion above Lemma \ref{L:ker_Jhat}, but for
    brevity we use this compact symbolic notation here.} $\bJ := D_\bA
  \bF(\bA)$, where $F_j(\bA):=\Omega A_j -\cN_j(\bA)$, has a simple
  zero eigenvalue.
\end{definition}
Note that due to the phase invariance $\bA\mapsto e^{\ri \nu}\bA, \nu
\in \R$ of \eqref{E:CME_NLB} the Jacobian is singular.
\begin{definition}\label{D:reverse}
  $\bA \in \C^N$ is \textit{reversible} if
  \[\bA\in \Vrev = \{\bv\in \C^N : v_i= \overline{v}_{i'} \text{ for
    all } i\in\{1,\dots,N\}\},\] where $i'$ is given by
  $\B\ni k^{(i')} \dot{=} -k^{(i)}$. 
\end{definition}
Reversibility is a symmetry of the solution. The
motivation for restricting to reversible non-degenerate solutions
$\bA$ is to ensure the invertibility of $\bJ$ in the fixed point
iteration for the singular part of the Lyapunov-Schmidt decomposition
in \S\ref{S:sing}. Within $\Vrev$ the phase invariance is, indeed, no
longer present. The choice of $\Vrev$ in Definition \ref{D:reverse} is
natural and based on the intrinsic symmetry \eqref{E:sym_Bloch} of the
Bloch eigenfunctions which ensures the $j\mapsto j'$ complex
conjugation symmetry among the coefficients in \eqref{E:CME_NLB} and,
hence, the possibility of reversible solutions. 
Note 
that \eqref{E:sym_Bloch} follows directly from $V(x)\in\R$. 

Next, we assume (H1-H6) and use the Lyapunov-Schmidt decomposition in Bloch
variables together with the Banach fixed point theorem to prove the main result, i.e., Theorem \ref{T:main}, which justifies the
formal asymptotics for solutions at $\omega=\omega_*+\Omega \eps^2$.

\subsection{Lyapunov-Schmidt Decomposition} \label{S:LS}

Due to the completeness of the Bloch waves $(\xi_n(\cdot,k))_{n\in
  \N}$ in $L^2(\P)$ we can expand \beq\label{E:phi_j}
\phi_j(x)=\sum_{n\in\N}\Phi^{(j)}_{n} \xi_n(x,k^{(j)}) \ \text{ with }
\ \Phi_n^{(j)} = (\phi_j(\cdot),\xi_n(\cdot,k^{(j)}))_{L^2(\P)}\in
\C.  \eeq As the following lemma shows, working with $\phi_j$ in the
$H^s(\P)$ space is equivalent to working with
$\Phibf^{(j)}:=(\Phi^{(j)}_n)_{n\in \N}\in l_{s/d}^2$, where
$$l^2_{s/d}=\{{\bf F}=(F_n)_{n\in \N} \in l^2: \|{\bf F}\|_{l^2_{s/d}}^2=\sum_{n\in \N}(1+n)^{2s/d}|F_n|^2<\infty\}.$$
\blem\label{L:norm_equi_Hs} For $s\geq 0$ the following norm
equivalence holds.  There exist constants $C_1,C_2>0$ such that
\[C_1 \|f\|_{H^s(\P)} \leq \|{\bf F}\|_{l^2_{s/d}} \leq C_2
\|f\|_{H^s(\P)} \quad \text{for all } f\in H^s(\P),\] where ${\bf
  F}:=(F_n)_{n\in \N}$ is related to $ f \in H^s(\P)$ by
\eqref{E:phi_j}.  
\elem 
The proof is analogous to that of Lemma 3.3 in
\cite{BSTU06}, see also \cite[\S4.1]{DU09}. The main ingredients are
firstly the fact that for $c>0$ large enough (such that $c+\omega_n(k)
>0$ for all $n$ and $k$, e.g. $c> -\essinf V$) the squared norm
$\|f\|_{H^s(\P)}^2$ is equivalent to
$$\int_{\R^d}\left|(c-\Delta+V(x))^{s/2}f(x)\right|^2 dx=\sum_{n\in \N}(c+\omega_n(k))^s\|p_n(\cdot,k)\|^2_{L^2(\P)}|F_n|^2 = \sum_{n\in \N}(c+\omega_n(k))^s|F_n|^2.$$ 
Secondly, one uses the asymptotic distribution of bands $\omega_n(k)$
in $d$ dimensions, see \cite[p.55]{Hoermander3}: 
there are constants $c_1,c_2,c_3>0$ such that
\beq\label{E:as_bands}
c_1n^{2/d}\leq \omega_n(k)+c_3 \leq  c_2 n^{2/d} \quad \forall n\in \N \ \forall k\in \B.
\eeq

\medskip

For the subsequent analysis we define for each $k^{(j)}\in K$ the set $\cAt_j$ of indices producing
$k^{(j)}$ through the nonlinearity analogously to the definition of
$\cA_j$, i.e.
\[\cAt_j:=\{(\alpha, \beta, \gamma)\in \{1,\dots,M\}^3:
k^{(\alpha)}-k^{(\beta)}+k^{(\gamma)}- k^{(j)} \in\Z^d\}.\] 

For the ansatz \eqref{E:NLB_ans}, \eqref{E:phi_j} equation
\eqref{E:SNLS_dD} is equivalent to the algebraic system
\beq\label{E:SNLS_bloch}
\cF^{(j)}_n(\vec{\Phibf}):=(\omega_n(k^{(j)})-\omega_* - \Omega
\eps^2)\Phi^{(j)}_n + \sigma G^{(j)}_n =0, \quad j\in \{1,\dots,M\}, \
n\in \N, \eeq where
\begin{align*}
  G^{(j)}_n &= \langle g_j, \xi_n(\cdot, k^{(j)})\rangle_{L^2(\P)}=\int_{\P}g_j(x) \overline{\xi_n(x, k^{(j)})}dx,\\
  g_j(x) &= \sum_{(\alpha,\beta,\gamma)\in \cAt_j} \sum_{n,o,q\in \N}
  \Phi^{(\alpha)}_n \overline{\Phi^{(\beta)}_o}\Phi^{(\gamma)}_q
  \xi_n(x,k^{(\alpha)})
  \overline{\xi_o(x,k^{(\beta)})}\xi_q(x,k^{(\gamma)}) =
  \sum_{(\alpha,\beta,\gamma)\in \cAt_j}
  \phi_\alpha\overline{\phi_\beta}\phi_\gamma.
\end{align*}

Due to the kernel of the linear multiplication operator at $\eps=0$ in
\eqref{E:SNLS_bloch} we use a \textit{Lyapunov-Schmidt
  decomposition} in order to characterize the bifurcation from
$\omega=\omega_*$ (i.e. from $\eps=0$). For $j\in \{1,\dots,M\}$ we let 
$$\ds I(j):= \begin{cases}
      \N \setminus \{n_j\} \ &\mbox{if}\ 1\leq j\leq N \\ \N \
      &\mbox{if}\ j > N \end{cases}, \text{ and let }
I_R:= \left\{(j,n): j\in \{1,\dots,M\}, n\in I(j)\right\}$$ 
and write  
\begin{align*}
  \phi(x) &= \eps  \phising(x) + \psi(x), \qquad \phising(x) = \sum_{j=1}^N B_j \xi_{n_j}(x,k^{(j)}), \quad \psi(x) = \sum_{(j,n)\in I_R} \Psi_n^{(j)}\xi_n(x,k^{(j)})
\end{align*}
with $0<\eps \ll 1$, $B_j\in \C$ and $ \Psibf^{(j)}
:=(\Psi_n^{(j)})_{n\in \N} \in l_{s/d}^2$.  In other words we set
\beq\label{E:LS_bloch}
\Phibf^{(j)} = \begin{cases} \eps B_j e_{n_j} + \Psibf^{(j)} \quad &\text{with} \ \Psibf^{(j)}\in l_{s/d}^2, \Psi^{(j)}_{n_j}=0  \ \mbox{ for } 1\leq j\leq N\\
  \Psibf^{(j)} \quad &\text{with} \ \Psibf^{(j)}\in l_{s/d}^2 \ \mbox{
    for } j >N,
\end{cases}
\eeq where $e_{n_j}$ is the $n_j$-th Euclidean unit vector in
$\R^\N$. Analogously to $\phi_j$ we also define
\[\psi_j:= \sum_{n\in I(j)} \Psi_n^{(j)}\xi_n(x,k^{(j)}).\]

This decomposition splits problem \eqref{E:SNLS_bloch} into
\begin{align}
  \cF^{(j)}_n:=(\omega_n(k^{(j)})-\omega_* - \Omega \eps^2)\Psi^{(j)}_n +\sigma G^{(j)}_n &=0, \qquad (j,n)\in I_R,\label{E:LS_reg}\\
  \cF^{(j)}_{n_j}:=-\eps^3 \Omega B_j +\sigma G^{(j)}_{n_j}&=0, \qquad
  j\in \{1,\dots,N\} \label{E:LS_sing}.
\end{align}

The following program is analogous to that in \cite{DPS09,DU09}.
Namely, for $(B_1,\dots, B_N)\in \C^N$ given, we first show the existence of a
small solution $(\Psibf^{(j)})_{j\in \N}$ of the regular part
\eqref{E:LS_reg} and then prove a persistence result relating certain
(reversible and non-degenerate) solutions $(A_1, \dots, A_N)\in \C^N$
of \eqref{E:CME_NLB} to solutions $(B_1, \dots, B_N)\in \C^N$ of the
singular part \eqref{E:LS_sing} including an estimate on their
difference, and finally provide an estimate of $\|\phi-\eps
\sum_{j=1}^NA_j\xi_{n_j}(\cdot, k_*^{(j)})\|_{H^s(\P)}$.

\section{Regular Part of the Lyapunov-Schmidt
  Decomposition}\label{S:reg}
We define the following spaces and norms
\begin{align*}
  \cS(s) := &\biggl\{\phi = \sum_{j\in \N} \phi_j:  \phi_j\in H^s(\P) \ \forall j, \|\phi\|_{\cS(s)}:= \sum_{j\in \N} \|\phi_j\|_{H^s(\P)} < \infty \ \text{ and } \\
  & \forall j \exists k \in \B \text{ such that } \phi_j(x+2\pi e_m) = e^{\ri 2\pi k_m}\phi_j(x), m=1,\dots,d \ \text{for a.e.} \ x\in \R^d \biggr\}\\
  \cX(s) := &\biggl\{\vec{\Phibf} = (\Phibf^{(j)})_{j\in \N}: \quad
    \|\vec{\Phibf}\|_{\cX(s)}:=\sum_{j\in\N}
    \|\Phibf^{(j)}\|_{l^2_{s/d}}<\infty\biggr\}.
\end{align*}

Note that the condition $k\in \B$ in the definition of $\cS(s)$ can be
replaced by $k\in \R^d$ because each $k\in \R^d$ can be written as
$k=\tilde{k}+\kappa$, where $\tilde{k}\in \B$ and $\kappa\in
\Z^d$. Also note that $\vec{\Phibf}$ is a sequence of sequences. 
Similarly we denote
$$\vec{\Psibf}:=(\Psibf^{(j)})_{j=1}^M \ \text{and} \  \vec{\Gbf} := (\Gbf^{(j)})_{j=1}^M.$$ 
Clearly, the ansatz \eqref{E:NLB_ans} satisfies $\phi \in \cS(s)$ if
and only if $\phi_j\in H^s(\P)$ for all $j\in\{1,\dots,M\}$. Therefore, for the problem at hand, where the solution consists of $M<\infty$ components $\phi_j$, the spaces $\cS(s)$ and $\cX(s)$ could be defined with finite sums over $j$. However, since the use of infinite sums in the definitions does not increase the complexity and since it may prove useful in future work on the case of irrational coordinates of $k_*^{(j)}$, we keep these general definitions.

We will need the following two lemmas, the first following 
directly from Lemma \ref{L:norm_equi_Hs}.
\blem\label{L:norm_equi} For
$s\geq 0$ there exist $c_1,c_2>0$ such that for all 
$$
\cS(s)\ni \phi(\cdot) = \sum_{j\in \N}\sum_{n\in \N} \Phi_n^{(j)}
\xi_n(\cdot,k^{(j)})
\text{ we have }
c_1 \|\phi\|_{\cS(s)} \leq \|\vec{\Phibf}\|_{\cX(s)} \leq c_2
\|\phi\|_{\cS(s)}.
$$
\elem 
\blem\label{L:algebra} For $s>d/2$ the space $\cS(s)$ is an algebra,
i.e. there is a constant $c>0$ such that $\|fg\|_{\cS(s)}\leq
c\|f\|_{\cS(s)}\|g\|_{\cS(s)}$ for all $f,g\in \cS(s)$.  \elem
\begin{proof}
  We define the sets $K_f$ and $K_g$ of $k-$points, which determine
  the quasiperiodicity of the functions $f_j$ and $g_j, j\in \N$, i.e.
  \begin{align*}
    K_f &:= \{k\in \B: \ \exists j\in\N \mbox{ with } f_j(x+2\pi e_m) =e^{2\pi\ri k_m}f_j(x) \text{ for all } m=1,\dots,d \text{ and a.e. } x\in  \R^d\},\\
    K_g &:= \{k\in \B: \ \exists j\in\N \mbox{ with } g_j(x+2\pi e_m)
    =e^{2\pi\ri k_m}g_j(x) \text{ for all } m=1,\dots,d \text{ and a.e. } x\in \R^d
    \}.
  \end{align*}
  We have
  \begin{align*}
    \|fg\|_{\cS(s)} &=\left\|\left(\sum_{\alpha\in \N} f_\alpha\right)\left(\sum_{\beta\in \N} g_\beta \right)\right\|_{\cS(s)} = \sum_{\stackrel{k^{(j)}\in K_f+ K_g}{k^{(j)} \text{ distinct}}} \left\|\sum_{k^{(\alpha)}+k^{(\beta)}\in k^{(j)}+ \Z^d} f_\alpha g_\beta \right\|_{H^s(\P)} \\
    &\leq c \sum_{\stackrel{k^{(j)}\in K_f+ K_g}{k^{(j)} \text{
          distinct}}} \sum_{k^{(\alpha)}+k^{(\beta)}\in k^{(j)} +
      \Z^d} \|f_\alpha\|_{H^s(\P)}\|g_\beta \|_{H^s(\P)} =
    c\|f\|_{\cS(s)}\|g\|_{\cS(s)},
  \end{align*}
  where the inequality follows by the triangle inequality and by the
  algebra property of the $H^s$ norm
$$
\|uv\|_{H^s(\P)}\leq C \|u\|_{H^s(\P)}\|v\|_{H^s(\P)} \quad
\forall u,v\in H^s(\P) \text{ and } s>d/2,
$$
see Theorem 5.23 in \cite{Adams}.
\end{proof}

Our result on the regular part of the Lyapunov-Schmidt decomposition
is the following 
\bprop\label{P:reg} 
Assume (H1)-(H5) and let
$\bB:=(B_1,\dots,B_N)\in \C^N$ be given (not necessarily a solution of
\eqref{E:LS_sing}). There exist $\eps_0>0$ and $C=C(\|\bB\|_{l_1})>0$ such for all $\eps\in(0,\eps_0)$ 
there exists a
solution $\vec{\Psibf} \in \cX(s)$ of \eqref{E:LS_reg} such that
\[\|\vec{\Psibf}\|_{\cX(s)} \leq C \eps^3.\] 
\eprop 
\bpf Writing \eqref{E:LS_reg} in the fixed point
formulation
\[\Psi^{(j)}_n = (\omega_n(k^{(j)})-\omega_*)^{-1}(\eps^2\Omega
\Psi^{(j)}_n - \sigma G^{(j)}_n)=:H^{(j)}_n(\vec{\Psibf}), \quad
(j,n)\in I_R,\] we seek a fixed point with $\|\vec{\Psibf}\|_{\cX(s)}
\leq \text{const.} \eps^3$.  Lemma \ref{L:norm_equi} allows us to work
interchangeably in $\cS(s)$ in the physical variables.  We show the
contraction property of $\vec{\bf H}$ within
\[D_{C \eps^3}:=\{\vec{\Psibf}: \ \|\vec{\Psibf}\|_{\cX(s)}\leq
C\eps^3\}\] for some $C>0$.

The nonlinearity is
\[|\phi|^2\phi = \eps^3|\phising|^2\phising +\eps^2
\left(2|\phising|^2\psi +\phising^2\overline{\psi}\right) +\eps
\left(2\phising|\psi|^2+\overline{\phising}\psi^2\right)+|\psi|^2\psi\]
such that we need to bound terms of the form
$\eps^3|\phising|^2\phising$, $\eps^2|\phising|^2\psi$, $\eps
\phising|\psi|^2$, and $|\psi|^2\psi$. Using the algebra property from
Lemma \ref{L:algebra} and the regularity of Bloch waves, we obtain 
\begin{align*}
  \eps^3\||\phising|^2\phising\|_{\cS(s)} &\leq c\eps^3
  \|\phising\|_{\cS(s)}^3 \leq c\eps^3 \left(\sum_{j=1}^N|B_\alpha|
    \|\xi_{n_j}(\cdot,k^{(j)})\|_{H^s(\P)}\right)^3 \leq c \eps^3
  \|\bB\|_{l^1}^3.
\end{align*} 
Similarly, for the remaining terms we have 
\begin{align*}
  \eps^2 \| |\phising|^2\psi\|_{\cS(s)} & \leq c\eps^2 \|\bB\|_{l^1}^2 \|\psi\|_{\cS(s)},\\
  \eps\|\phising |\psi|^2\|_{\cS(s)} & \leq c \eps \|\bB\|_{l^1} \|\psi\|^2_{\cS(s)},\\
  \||\psi|^2\psi\|_{\cS(s)} &\leq c\|\psi\|^3_{\cS(s)}.
\end{align*}

Next, thanks to assumptions (H3)-(H5) we have the uniform lower bound
$$|\omega_n(k^{(j)})-\omega_*|>c>0 \ \text{for all} \ (j,n)\in I_R.$$
>From (H4) follows that the $j-$set in $I_R$ is finite so that
the minimum of $|\omega_n(k^{(j)})-\omega_*|$ in $j$ can be
taken. (H3) and (H5) ensure that the minimum is positive.

Collecting the above estimates, we thus have
\[\|\vec{\bf H}\|_{\cX(s)} \leq C\left[\eps^3\|\bB\|_{l^1}^3 +\eps^2
  (\|\bB\|_{l^1}^2+|\Omega|) \|\vec{\Psibf}\|_{\cX(s)}+\eps
  \|\bB\|_{l^1} \|\vec{\Psibf}\|_{\cX(s)}^2
  +\|\vec{\Psibf}\|_{\cX(s)}^3\right].\] We conclude that for $\eps>0$
small enough $\vec{\bf H}$ maps $D_{C\eps^3}$ to itself.

Similarly, the contraction property of ${\bf H}$ follows by the same
estimates as above, the simple identities
\begin{align*}
  |\psi_a|^2-|\psi_b|^2 &= \tfrac{1}{2}\left[(\psi_a-\psi_b)(\overline{\psi_a}+\overline{\psi_b})+(\psi_a+\psi_b)(\overline{\psi_a}-\overline{\psi_b})\right],\\
  |\psi_a^2-\psi_b^2|&= |\psi_a+\psi_b||\psi_a-\psi_b|,\\
  |\psi_a|^2\psi_a-|\psi_b|^2\psi_b &=
  (|\psi_a|^2+|\psi_b|^2)(\psi_a-\psi_b)+\psi_a\psi_b(\overline{\psi_a}-\overline{\psi_b}),
\end{align*}
and by the algebra property. We find
\begin{align*}
  \|\vec{\bf H}(\vec{\Psibf}_a)-\vec{\bf H}(\vec{\Psibf}_b)\|_{\cX(s)} \leq & C\left[\eps^2 (\|\bB\|_{l^1}^2+|\Omega|) +\eps \|\bB\|_{l^1} (\|\vec{\Psibf}_a\|_{\cX(s)}+\|\vec{\Psibf}_b\|_{\cX(s)}) \right.\\
  &\left. +\|\vec{\Psibf}_a\|^2_{\cX(s)}+\|\vec{\Psibf}_b\|^2_{\cX(s)}\right]\|\vec{\Psibf}_a-\vec{\Psibf}_b\|_{\cX(s)}
\end{align*}
for all $\vec{\Psibf}_a,\vec{\Psibf}_b\in \cX(s)$.  In conclusion, the
existence of a solution $\vec{\Psibf} \in D_{C(\|\bB\|_{l_1})\eps^3}$ follows.
\epf

\section{Singular Part of the Lyapunov-Schmidt Decomposition,
  Persistence}\label{S:sing}
The singular part \eqref{E:LS_sing} of the Laypunov-Schmidt
decomposition is equivalent to the \textit{extended algebraic coupled
  mode equations} \beq\label{E:ECME} \Omega B_j - \cN_j(B_1, \dots,
B_N) = R_j, \qquad j\in \{1,\dots,N\} \eeq with $R_j:=\eps^{-3}
G^{(j)}_{n_j}-\cN_j(B_1,\dots, B_N)$. Proposition \ref{P:reg} thus leads to the following
\bcor\label{C:ECME} Assume(H1)-(H5), and let 
$(B_1,\dots,B_N)\in \C^N$ be a solution of \eqref{E:ECME}.  
There exist $\eps_0>0$ and $C>0$ such that for all $\eps\in(0,\eps_0)$ 
equation \eqref{E:SNLS_dD} with $\omega=\omega_*+\eps^2\Omega$ has a
nonlinear Bloch wave solution $\phi$ of the form \eqref{E:NLB_ans}
such that
\[\left\|\phi(\cdot)-\eps\sum_{j=1}^N
  B_j\xi_{n_j}(\cdot,k^{(j)})\right\|_{H^s(\P)}\leq C \eps^3. 
\]
\ecor 
Corollary \ref{C:ECME} is of
little practical use since $G_{n_j}^{(j)}$ in \eqref{E:ECME} depend on the unknown $\psi$ such that solving \eqref{E:ECME} for
$(B_1, \dots, B_N)$ explicitly is not possible. This problem can be
avoided by showing persistence of solutions $(A_1, \dots, A_N)\in
\C^N$ of the formally derived explicit ACMEs \eqref{E:CME_NLB} to
solutions $(B_1,\dots,B_N)\in \C^N$ of \eqref{E:ECME}, which is our
next step. We show that persistence holds for ``reversible
non-degenerate'' solutions $(A_1,\dots,A_N)\in \C^N$. The problem then
reduces to finding reversible non-degenerate solutions of the ACMEs.
Writing the ACMEs as $F_j(A_1,\dots,A_N)=0$, equation
\eqref{E:ECME} reads
\[F_j(B_1,\dots,B_N) = R_j, \qquad j\in \{1,\dots,N\}.\]
\blem\label{L:R_est} Assume (H1). Given $\vec{\Psibf}\in \cX(s)$ with $\|\vec{\Psibf}\|_{\cX(s)}<C \eps^3$ we have
\[|R_j| \leq C \eps^2\] for all $j\in \{1,\dots,N\}$, where
$C=C(\|\bB\|_{l^1})>0.$ \elem 
\bpf 
Substituting for $\phi_j$, the
decomposition \eqref{E:LS_bloch}, we get
\begin{align*}
  R_j = & ~\eps^{-3}\sigma G^{(j)}_{n_j} - \cN_j(B_1,\dots, B_N) =\\
  = & ~ \eps^{-1}\sigma \left\{ 2\sum_{(\alpha,\beta,\gamma)\in \cA_j} B_\alpha\overline{B}_\beta\langle\psi_\gamma(\cdot)\xi_\alpha(\cdot,k^{(\alpha)})\overline{\xi_\beta(\cdot,k^{(\beta)})},\xi_{n_j}(\cdot,k^{(j)})\rangle_{L^2(\P)}\right.\\
  & \qquad \left. +\sum_{(\alpha,\beta,\gamma)\in \cA_j} B_\alpha
    B_\gamma\langle\overline{\psi}_\beta(\cdot)\xi_\alpha(\cdot,k^{(\alpha)})\xi_\gamma(\cdot,k^{(\gamma)}),\xi_{n_j}(\cdot,k^{(j)})\rangle_{L^2(\P)}\right\}
\end{align*}
\begin{align*}
  \qquad &+\eps^{-2}\sigma \left\{ 2\sum_{(\alpha,\beta,\gamma)\in \cA_j} B_\alpha\langle\overline{\psi}_\beta(\cdot)\psi_\gamma(\cdot)\xi_\alpha(\cdot,k^{(\alpha)}),\xi_{n_j}(\cdot,k^{(j)})\rangle_{L^2(\P)}\right.\\
  & \qquad \left. +\sum_{(\alpha,\beta,\gamma)\in \cA_j} \overline{B}_\beta \langle\psi_\alpha(\cdot)\psi_\gamma(\cdot)\overline{\xi_\beta(\cdot,k^{(\beta)})},\xi_{n_j}(\cdot,k^{(j)})\rangle_{L^2(\P)}\right\}\\
  & + \eps^{-3}\sigma \sum_{(\alpha,\beta,\gamma)\in \cA_j}
  \langle\psi_\alpha(\cdot)\overline{\psi}_\beta(\cdot)\psi_\gamma(\cdot),\xi_{n_j}(\cdot,k^{(j)})\rangle_{L^2(\P)}.
\end{align*}
With the Cauchy-Schwarz inequality, the regularity of Bloch waves, and
using the estimate $\|\psi_\alpha\|_{H^s(\P)}\leq C \eps^3$ for all $\alpha$, which follows from the assumption $\|\vec{\Psibf}\|_{\cX(s)}<C \eps^3$, we obtain the desired estimate for $|R_j|$.
\epf

Next we let $\bB =\bA+\bb$, where similarly to $\bA$ we denote $\bB:=(B_1,\dots, B_N)^T$. The difference $\bb$ solves
\beq \label{E:b_eq} \bJ\bb =\bW(\bb), \qquad \bW(\bb) :=
\bR(\bA+\bb)-(\bF(\bA+\bb)-\bJ\bb), \eeq where $\bF:=(F_1,\dots,
F_N)^T$, $\bR:=(R_1,\dots, R_N)^T$, and $\bJ = D_\bA \bF(\bA)$ is the Jacobian\footnote{A symbolic notation for the Jacobian used again.} of $\bF$ at $\bA$. 
Due to $\bF(\bA)=0$, we get
that $\bF(\bA+\bb)-\bJ\bb$ is at least quadratic in $\bb$ so that for
$|\bb|$ small we have (in the Euclidean norm $|\cdot|$)
\[|\bF(\bA+\bb)-\bJ\bb|\leq c |\bb|^2.\] As a result
\beq\label{E:W_est} |\bW(\bb)|\leq
c\left\{\eps^2+\eps^2|\bb|+|\bb|^2\right\} \eeq for $|\bb|$ small,
where the $c\eps^2$ term comes from $\bA$-homogenous terms in $\bR$
and $\eps^2|\bb|$ from linear terms in $\bb$.  

We aim to apply a fixed point 
argument
on $\bb =\bJ^{-1}\bW(\bb)$ in a neighborhood of $0$ to 
produce a solution $\bb$ with $|\bb|<c\eps^2$. However, due to the phase
invariance $\bA\mapsto e^{\ri \nu} \bA, \nu\in \R$ of $\bF(\bA)=0$ the
Jacobian $\bJ$ is not invertible.
To overcome this difficulty, we assume the non-degeneracy of $\bA$, see Definition \ref{D:nondegen}.
Second, we restrict
$\bA$ and $\bb$ to the reversible space $\Vrev$, see Definition \ref{D:reverse}, in which
$\bJ$ is invertible, as shown below. Our precise requirements on $\Vrev$ are:
\beq\label{E:rev_cond}
\text{If } 0\neq \bA\in \Vrev, \text{ then} \begin{cases} \text{(i)} \ &\exists \delta>0 \text{ such that }  |\bJ \bb| > \delta |\bb| \text{ for all } \bb \in \Vrev, \\
  \text{(ii)} \ &\bJ^{-1}\bW(\bb)\in \Vrev \text{ for all } \bb \in
  \Vrev.\end{cases} 
	\eeq 
	
To check (i) and (ii) in \eqref{E:rev_cond}, we first formulate
$\bb,\bA,\bF$ and $\bJ$ in real variables and define the symmetry matrix
$\hat{S}$ corresponding to the reversibility symmetry in $\Vrev$. For
$\bv\in \C^N$ define $\hat{\bv}:=\bspm\bv_R \\ \bv_I\espm \in
\R^{2N}$, where $\bv_R \in \R^N$ and $\bv_I\in \R^N$ are the vectors
of real and imaginary parts of $\bv$. Then \beq\label{E:Vrev_S} \bv\in
\Vrev \quad \Leftrightarrow \quad \hat{\bv} = \hat{S} \hat{\bv}, \eeq
where
$$\hat{S}=\begin{pmatrix}P & \\ & -P\end{pmatrix}, \quad P = (\be_{1'},\be_{2'},\dots, \be_{N'}),
$$
and $\be_i$ is the $i-$th Euclidean unit vector in $\R^N$.
Let us denote by $\hat{\bA},\hat{\bb},\hat{\bF} \in \R^{2N}$ the
quantities $\bA,\bb,\bF$ in real variables and let $\hat{\bJ}\in
\R^{2N\times 2N} = D\hat{\bF}$ be the Jacobian of $\hat{\bF}$.

The uniform boundedness property (i) in \eqref{E:rev_cond} follows
since by the non-degeneracy condition $\hat{\bJ}$ has only one zero
eigenvalue and for $\bb\in\Vrev$ is $\hat{\bb}$ orthogonal to the
corresponding eigenvector. This is shown in the following 
\blem\label{L:ker_Jhat} If
$\bb,\bA\in \Vrev, \bF(\bA)=0$, and if $\bA$ is non-degenerate, then
$$\hat{\bb}^T \hat{\bv} =0 \ \text{for all} \ \hat{\bv} \in \ker(\hat{\bJ}) 
= \text{span}\left\{\bspm 0 & -I\\ I & 0 \espm \hat{\bA} \right\}.$$
\elem 
\bpf The well known fact $\ker(\hat{\bJ}) =
\text{span}\left\{\bspm 0 & -I\\ I & 0 \espm \hat{\bA} \right\}$
follows from the phase invariance $\bF(e^{\ri \nu}\bA)=0$ for all $\nu
\in \R$ by rewriting it in real variables, differentiating in $\nu$
and evaluating at $\nu=0$.  Using now \eqref{E:Vrev_S} for $\hat{\bA}$
and $\hat{\bb}$, we get
\[\hat{\bb}^T \bpm 0 & -I\\ I & 0 \epm \hat{\bA} = \hat{\bb}^T \bpm
P^T & 0\\ 0 & -P^T \epm \bpm 0 & P\\ P & 0 \epm \hat{\bA} = -
\hat{\bb}^T \bpm 0 & -I\\ I & 0 \epm \hat{\bA}.\] 

\vspace{-8mm}
\epf

For (ii) in \eqref{E:rev_cond} let us first show that $\bA,\bb \in
\Vrev \ \Rightarrow \ \bW(\bb)\in \Vrev$. Because of the symmetry
\eqref{E:sym_Bloch} and the
symmetry $(\alpha,\beta,\gamma)\in \cA_j\Leftrightarrow
(\alpha',\beta',\gamma')\in \cA_{j'}$ we get from \eqref{E:mu} that
$\mu_{\alpha',\beta',\gamma',j'}=\overline{\mu_{\alpha,\beta,\gamma,j}}$
for all $\alpha,\beta,\gamma,j\in\{1,\dots,N\}.$ As a result, $\bF$ has
the symmetry \beq\label{E:Fsym}
\hat{\bF}(\hat{S}\hat{\bv})=\hat{S}\hat{\bF}(\hat{\bv}) \quad
\mbox{for all } \hat{\bv}\in \R^{2N}.  \eeq For $\bA,\bb\in \Vrev$ 
this results in
$\hat{\bF}(\hat{\bA}+\hat{\bb}) =
\hat{S}\hat{\bF}(\hat{\bA}+\hat{\bb})$, i.e.~$\bF(\bA+\bb)\in \Vrev$. 

Next, differentiating \eqref{E:Fsym}, we get
\[\hat{\bF}'(\hat{S}\hat{\bv})\hat{S}=\hat{S}\hat{\bF}'(\hat{\bv})
\quad \mbox{for all } \hat{\bv}\in \R^{2N}.\] If $\bA\in \Vrev$, then
this translates for $\bv=\bA$ to \beq\label{E:Jsym}
\hat{\bJ}\hat{S}=\hat{S}\hat{\bJ} \eeq and for $\bA,\bb\in\Vrev$ we
thus have
$\hat{\bJ}\hat{\bb}=\hat{\bJ}\hat{S}\hat{\bb}=\hat{S}\hat{\bJ}\hat{\bb}$,
so that $\bJ\bb\in \Vrev$.

The last term in $\bW$ is $\bR$, where $R_j =
\eps^{-3}G^{(j)}_{n_j}-\cN_j$, $j=1,\dots,N$. For $\cN_j$ the above
identity
$\mu_{\alpha',\beta',\gamma',j'}=\overline{\mu_{\alpha,\beta,\gamma,j}}$
implies that for $\bA,\bb\in \Vrev$
$$\cN_j=\overline{\cN_{j'}}.$$

For $G^{(j)}_{n_j}$ we argue as follows. First, we define the symmetry
map
$$S:\vec{\Phibf}\mapsto S\vec{\Phibf}, \text { where } (S\vec{\Phibf})^{(j)} = \overline{\Phibf^{(j')}}.$$
\blem\label{L:Gsym} $\vec{\Gbf}$ commutes with $S$, i.e.
$$\vec{\Gbf}(S\vec{\Phibf}) = S\vec{\Gbf}(\vec{\Phibf}).$$
\elem \bpf For all $j\in \{1,\dots,M\}$ and $n\in \N$ we get, using
\eqref{E:sym_Bloch},
\begin{align*}
  G_{n}^{(j)}(S\vec{\Phibf}) &= \sum_{(\alpha\beta\gamma)\in \tilde{\cA}_j}\sum_{m,o,q\in \N}\overline{\Phi^{(\alpha')}_m}\Phi_o^{(\beta')}\overline{\Phi^{(\gamma')}_q}\int_{\P}\xi_m(x,k^{(\alpha)})\overline{\xi_o(x,k^{(\beta)})}\xi_q(x,k^{(\gamma)})\overline{\xi_n(x,k^{(j)})}dx\\
  &= \sum_{(\alpha'\beta'\gamma')\in \tilde{\cA}_{j'}}\sum_{m,o,q\in \N}\overline{\Phi^{(\alpha')}_m}\Phi_o^{(\beta')}\overline{\Phi^{(\gamma')}_q}\int_{\P}\overline{\xi_m(x,k^{(\alpha')})}\xi_o(x,k^{(\beta')})\overline{\xi_q(x,k^{(\gamma')})}\xi_n(x,k^{(j')})dx\\
  &= \overline{G_{n}^{(j')}(\vec{\Phibf}) }.
\end{align*}

\vspace{-8mm}
\epf 
\blem\label{L:psisym} If $\bB\in \Vrev,$ then there exists a
solution $\vec{\Psibf}$ of \eqref{E:LS_reg} with the properties as in
Proposition \ref{P:reg}, and such that
$$\vec{\Psibf} = S\vec{\Psibf}.$$
\elem 
\bpf Defining $\vec{\Phibf}_\text{sing}$ via
$$\vec{\Phibf}^{(j)}_\text{sing} = \begin{cases} B_j e_{n_j}, 
&j\in \{1,\dots,N\}, \\ 0, &j\in \{N+1,\dots, M\},\end{cases}$$
we have $\vec{\Phibf} = \eps \vec{\Phibf}_\text{sing}
+\vec{\Psibf}$. Due to (H6) is $\bB\in \Vrev$ equivalent to
$S\vec{\Phibf}_\text{sing} = \vec{\Phibf}_\text{sing}$.  And if
$S\vec{\Phibf}_\text{sing} = \vec{\Phibf}_\text{sing}$, then the fixed
point iteration $\vec{\Psibf} = \vec{\bf H}(\vec{\Psibf})$ preserves
the symmetry of $\vec{\Psibf}$, i.e.
$$\vec{\Psibf} = S\vec{\Psibf} \quad \Rightarrow \quad \vec{\bf H}(\vec{\Psibf}) = S \vec{\bf H}(\vec{\Psibf}).$$
This is clear from the form
$$H_n^{(j)}=(\omega_n(k^{(j)})-\omega_*)^{-1}\left(\eps^2\Omega \Psi_n^{(j)}-\sigma G_n^{(j)}(\eps \vec{\Phibf}_\text{sing} + \vec{\Psibf})\right)$$
and from Lemma \ref{L:Gsym}.  \epf

>From Lemma \ref{L:psisym} we conclude that given $\bB\in \Vrev$, the
full vector $\vec{\Phibf}$ is $S-$symmetric, i.e. $\vec{\Phibf} = \eps
\vec{\Phibf}_\text{sing} +\vec{\Psibf} = \eps
S\vec{\Phibf}_\text{sing} +S\vec{\Psibf}= S \vec{\Phibf}$. Lemma
\ref{L:Gsym} then yields for all $j\in \{1, \dots, N\}$
$$G^{(j)}_{n_j}=\overline{G^{(j')}_{n_j}}.$$
Thanks to (H6) $j'\in \{1,\dots,N\}$, and in conclusion $\bR\in \Vrev$
for $\bB=\bA+\bb\in \Vrev$.

Summarizing, we have $\bW\in \Vrev$ for $\bA,\bb\in\Vrev$. To conclude
the proof of (ii) in \eqref{E:rev_cond} we need to prove $\bv\in\Vrev
\Rightarrow \bJ^{-1}\bv\in \Vrev$. From \eqref{E:Jsym} we get within
$\Vrev$, where $\bJ^{-1}$ is defined,
\[\hat{S}\hat{\bJ}^{-1} \hat{S}^{-1} =\hat{\bJ}^{-1}.\]
If $\bv\in\Vrev$, then $\hat{\bv}=\hat{S}\hat{\bv}$ and
\[\hat{\bJ}^{-1} \hat{\bv} = \hat{S}\hat{\bJ}^{-1} \hat{S}^{-1}
\hat{S}\hat{\bv} = \hat{S}\hat{\bJ}^{-1} \hat{\bv}.
\] 
This shows that
$\bJ^{-1}\bv\in \Vrev$. We can thus finally solve the fixed point problem 
\eqref{E:b_eq} to obtain $\bb$ with $|\bb|< C \eps^2$. Herewith we obtain the following
\bprop\label{P:persist} Assume (H6) and let 
$\bA$ be a reversible non-degenerate solution of the coupled mode 
equations \eqref{E:CME_NLB}. There exist $\eps_0>0$ and $C>0$ such that 
for all $\eps\in(0,\eps_0)$ the following holds. 
Given $\vec{\Psibf}\in
\cX(s)$ with $\|\vec{\Psibf}\|_{\cX(s)} \leq C\eps^3$, 
there exists a
solution $\bB\in \Vrev$ of the extended coupled mode equations
\eqref{E:ECME} such that
\[|\bA-\bB| < C\eps^2.\] 
\eprop 
Our main result, i.e.~Theorem \ref{T:main}, for the bifurcation of
nonlinear Bloch waves 
follows from Corollary \ref{C:ECME}, Proposition \ref{P:persist} and
the triangle inequality.

\section{ACMEs for $N=1$ and $N=2$} \label{S:ACMEs_N12} 
We present here the complete solution structure of the ACMEs for the cases $N=1$ and
$N=2$.

\subsection{One Mode: $N=1$}
If $N=1$, then necessarily also $M=1$ since $S_3(\{k_*\})=\{k_*\}$ for
each $k_*\in \B$. Hence, $N=1$ is always consistent. However, only for
$k_*\in \{0,\tfrac{1}{2}\}^d$ condition (H6) is satisfied. The ACMEs
\eqref{E:CME_NLB} now have the scalar form \beq\label{E:CME_scalar}
\Omega A -\sigma \mu |A|^2 A=0, \quad \mu = \|\xi_{n_*}(\cdot,
k_*)\|_{L^4(\P)}^4>0, \eeq where $\xi_{n_*}(x,k_*)$ is the linear
Bloch wave for the selected eigenvalue index $n_*$. Note that $n_*$ has to
be chosen such that (H3) holds. Clearly, nonzero solutions of
\eqref{E:CME_scalar} satisfy
$$|A|=\sqrt{\tfrac{\Omega}{\sigma \mu}},$$
which implies a bifurcation to the left in
$\omega$ from $\omega_*$ in the focusing case $\sigma<0$ and to the
right in the defocusing case $\sigma >0$.

\subsection{Two Modes: $N=2$} \label{S:N2} Also for $N=2$ the solutions of the resulting
ACMEs can be calculated explicitly.  We discuss only solutions with
$A_1A_2\neq 0$. This is without any loss of generality because
if $k_*^{(2)}\in -k_*^{(1)} + \Z^d$, then the reversibility ${\bf
  A}\in \Vrev$ implies $A_2=\overline{A_1}$ and if $k_*^{(2)}\notin
-k_*^{(1)} + \Z^d$, then considering only one nonzero component
in ${\bf A}$ is equivalent to considering the case $N=1$.

For $N=2$ the form of the ACMEs depends on the choice of
$\{k_*^{(1)},k_*^{(2)}\}$. There are the following two cases.
\begin{enumerate}
\item[(a)] Let  
\beq\label{E:N2_case1} 2k_*^{(1)}-k_*^{(2)} \in
  k_*^{(2)} + \Z^d, \text{ i.e. } k_*^{(1)} \in k_*^{(2)}
  + \{-1/2,1/2\}^d.  \eeq This can be easily seen to be the
  consistent case $S_3(\{k_*^{(1)},k_*^{(2)}\}) \subset
  \{k_*^{(1)},k_*^{(2)}\} + \Z^d$, i.e. the case $M=N=2$. In this
  case we have
  \[\cA_1=\{(1,1,1),(1,2,2),(2,2,1),(2,1,2)\}, \
  \cA_2=\{(2,2,2),(2,1,1),(1,1,2),(1,2,1)\}, 
\] 
and the ACMEs read
\beq\label{E:N2_consist}
\begin{aligned}
    \Omega A_1 -\sigma\left[(\mu_{1111}|A_1|^2+2\mu_{1221}|A_2|^2)A_1+\mu_{2121}A_2^2\overline{A}_1\right]=&0,\\
    \Omega A_2
    -\sigma\left[(\mu_{2222}|A_2|^2+2\mu_{1221}|A_1|^2)A_2+\overline{\mu_{2121}}A_1^2\overline{A}_2\right]=&0,
  \end{aligned}\eeq
  where the obvious identities $\mu_{1221}=\mu_{2112}$ and
  $\mu_{1212}=\overline{\mu_{2121}}$ have been used. A simple
  calculation yields that solutions with both $A_1$ and $A_2$ nonzero
  satisfy
  \begin{align*}
    &\text{arg}(A_2)=\text{arg}(A_1)-\frac{\text{arg}(\mu_{2121})}{2}+q \frac{\pi}{2}, q\in \Z, \\
    &|A_1|^2 = \frac{\Omega}{\sigma}\frac{\gamma-\mu_{2222}}{\gamma^2
      -\mu_{1111}\mu_{2222}}, \ |A_2|^2 =
    \frac{\Omega}{\sigma}\frac{\gamma-\mu_{1111}}{\gamma^2
      -\mu_{1111}\mu_{2222}},
  \end{align*}
  where $\gamma := 2\mu_{1221}+(-1)^q|\mu_{2121}|$.

  A solution with $A_1,A_2\neq 0$ thus exists for
  $\sign(\Omega)=\sign(\sigma)$ if and only if
  $$\sign(\gamma-\mu_{2222})=\sign(\gamma-\mu_{1111})=\sign(\gamma^2
  -\mu_{1111}\mu_{2222})$$ is satisfied either for $q=0$ or $q=1$. For
  $\sign(\Omega)=-\sign(\sigma)$ the existence follows if and only if
  $$\sign(\gamma-\mu_{2222})=\sign(\gamma-\mu_{1111})=-\sign(\gamma^2
  -\mu_{1111}\mu_{2222})$$ either for $q=0$ or $q=1$.

  In order to satisfy the reversibility condition ${\bf A}\in \Vrev$,
  we need $A_2=\overline{A_1}$. This is possible if and only if
  $\mu_{1111}=\mu_{2222}$ such that $|A_1|=|A_2|$. The equality
  $A_2=\overline{A_1}$ then follows if we choose
$$\text{arg}(A_1)=\frac{\text{arg}(\mu_{2121})-q\pi}{4}.$$

\item[(b)] If \eqref{E:N2_case1} does
  not hold, then we have an inconsistent case $M>N=2$,
  \[\cA_1=\{(1,1,1),(1,2,2),(2,2,1)\}, \
  \cA_2=\{(2,2,2),(2,1,1),(1,1,2)\},\] 
and the ACMEs have the form
  \begin{align*}
    \Omega A_1 -\sigma (\mu_{1111}|A_1|^2+2\mu_{1221}|A_2|^2)A_1=&0,\\
    \Omega A_2 -\sigma (\mu_{2222}|A_2|^2+2\mu_{1221}|A_1|^2)A_2=&0.
  \end{align*}
  Solutions with both $A_1$ and $A_2$ nonzero satisfy
$$
|A_1|^2 =
\frac{\Omega}{\sigma}\frac{2\mu_{1221}-\mu_{2222}}{4\mu_{1221}^2
  -\mu_{1111}\mu_{2222}}, \ |A_2|^2 =
\frac{\Omega}{\sigma}\frac{2\mu_{1221}-\mu_{1111}}{4\mu_{1221}^2
  -\mu_{1111}\mu_{2222}}.$$ Again, the reversibility condition can be
satisfied (by choosing $\text{arg}(A_1)=-\text{arg}(A_2)$) if and only
if $\mu_{1111}=\mu_{2222}$.
\end{enumerate}

In one dimension $d=1$ with $N=2$ the only consistent cases satisfying
(H6) are
\[\{k_*^{(1)},k_*^{(2)}\} = \{0,1/2\} \text{ and }
\{k_*^{(1)},k_*^{(2)}\} = \{-1/4,1/4\}.\] In two dimensions $d=2$ with
$N=2$ there are 12 possible sets $\{k_*^{(1)},k_*^{(2)}\}$ 
satisfying (H6) and the consistency, namely 

\begin{tabular}{llll}
 $\left\{\bspm 0 \\[1mm]0 \espm, \bspm 1/2\\[1mm] 0 \espm\right\}$&
$\left\{\bspm 0 \\[1mm]0 \espm, \bspm 0\\[1mm]1/2\espm\right\}$&
$\left\{\bspm 0 \\[1mm]0 \espm, \bspm 1/2\\[1mm] 1/2 \espm\right\}$&
$\left\{\bspm 1/2 \\[1mm]0 \espm, \bspm 0\\[1mm] 1/2 \espm\right\}$\\[2mm]
$\left\{\bspm 1/2 \\0 \espm, \bspm 1/2\\ 1/2 \espm\right\}$&
$\left\{\bspm 0 \\1/2 \espm, \bspm 1/2\\ 1/2 \espm\right\}$&
$\left\{\bspm 1/4 \\ 0 \espm, \bspm -1/4\\ 0\espm\right\}$&
$\left\{\bspm 0\\ 1/4  \espm, \bspm 0 \\ -1/4\espm\right\}$\\[2mm]
$\left\{\bspm 1/4\\ 1/4  \espm, \bspm -1/4 \\ -1/4\espm\right\}$& 
$\left\{\bspm 1/4\\ -1/4  \espm, \bspm -1/4 \\ 1/4\espm\right\}$&
$\left\{\bspm 1/2 \\1/4 \espm, \bspm 1/2\\ -1/4 \espm\right\}$&
$  \left\{\bspm 1/4 \\1/2 \espm, \bspm -1/4\\ 1/2 \espm\right\}$.
\end{tabular}

\section{Numerical Examples in Two Dimensions $d=2$}
\label{S:numerics}
\def\pdep{{\tt pde2path}}

In the following numerical computations we use the package \texttt{pde2path}
\cite{p2p,p2p2} for numerical continuation and bifurcation 
in nonlinear elliptic systems of PDEs. The package uses linear finite
elements for the discretization, Newton's iteration for the
computation of nonlinear solutions and arclength 
continuation of solution branches. In the case $N=1$ below we
discretize $\P^2$ by $2*200^2=80000$ isosceles triangles of equal
size. For example B below with $N=2$ we use $2*280^2=156800$
triangles. This fine discretization is needed only in the tests of
$\eps$-convergence of the asymptotic error to ensure that the
asymptotic error dominates the discretization error. For all the numerical solutions (solution branches) presented in this and the following sections 
we verified that these approximate PDE solutions by standard 
error estimators and adaptive mesh--refinement.

For $N=1$ we simply write $\phi(x) =e^{\ri k_*\cdot x} \eta(x)$ and 
use real variables $\eta=u_1+\ri u_2$ to obtain 
\huga{\label{SGPb}
0=-\begin{pmatrix} \Delta u_1 \\ \Delta u_2\end{pmatrix}+2 \begin{pmatrix} k_*\cdot \nabla u_2\\ -k_*\cdot \nabla u_1\end{pmatrix} + (|k_*|^2-\omega +V(x)) \begin{pmatrix} u_1 \\ u_2 \end{pmatrix} +\sigma (u_1^2+u_2^2) \begin{pmatrix} u_1 \\ u_2 \end{pmatrix} 
}
on the torus $\mathbb{T}^2=\R^2/(2\pi\Z^2)$. For the consistent case 
with $N>1$ we may plug $\phi(x)=\sum_{j=1}^N e^{\ri k_*^{(j)}\cdot x}\eta_j(x)$ with $2\pi$--periodic $\eta_j$ 
into \reff{E:SNLS_dD} and collect 
terms multiplying $e^{\ri k_*^{(j)}\cdot x}$ in separate equations. Setting 
$\eta_j=u_1^{(j)}+\ri u_2^{(j)}$ we obtain a real system of 
$2N$ equations for $u=(u_1^{(1)}, u_2^{(1)},\ldots,u_1^{(N)},u_2^{(N)})$. 

We may then use two methods to generate branches of NLBs. 
The first is to let \pdep\ find the bifurcation points from 
the trivial branch $u=0$ and then perform branch switching 
to and continuation of the bifurcating branches. This is what we did 
in Example \ref{Ex:NLB_intro} from the Introduction to obtain 
Figure \ref{f2}. 
However, as  due to the phase invariance 
the eigenvalues of the linearization 
of \reff{SGPb} are always double, this needs 
some slight modification of the standard 
bifurcation detection and branch--switching 
routines of \pdep, see \cite[\S2.6.4]{p2p2b}. Thus,  in the examples below 
we alternatively use the asymptotic approximation $\phi(x)=\eps
\sum_{j=1}^N A_j\xi_{n_j}(x,k_*^{(j)})$ as the initial
guess in the Newton's iteration for the first continuation step near
$\omega=\omega_*$. 

We choose the potential \eqref{E:V}, which is the same  as in \cite{DU09}. 
The band structure along the boundary of the irreducible Brillouin
zone is plotted in Fig.~\ref{F:band_str}(b), 
and in Example \ref{Ex:NLB_intro} we already gave an overview of 
the lowest bifurcations at point $X$ with $N=1$. 
In the following examples we consider in more detail the 
points marked (A),(B),(C). Note that (C) is not a case of high symmetry points as $k_*^{(1)}=-k_*^{(2)}=(1/4,1/4)\notin \{0,1/2\}^2$. 

\subsection{Numerical Example for $N=1$.}\label{num1s}
For $N=1,d=2$ the only cases which satisfy (H6) are
$$k_*=(0,0), k_*=(1/2,0), k_*=(0,1/2), \text{ and } k_*=(1/2, 1/2).$$ 
$k_*=(1/2,0)$ with $N=1$ was considered in Example \ref{Ex:NLB_intro}, and in 
\S\ref{num2s} we reconsider this $k_+$ at point $(B)$ in Fig.~\ref{F:band_str} with $N=2$. 
Here we present in some more detail nonlinear Bloch waves bifurcating from 
point $(A)$ with $k_*=(1/2,1/2)$. 

\begin{figure}[!ht]
  \begin{center}
\scalebox{0.5}{\includegraphics{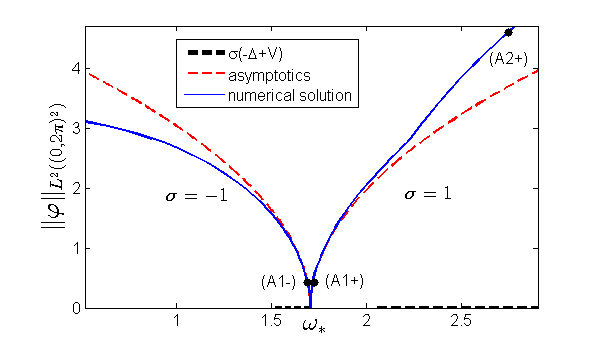}}
\scalebox{0.53}{\includegraphics{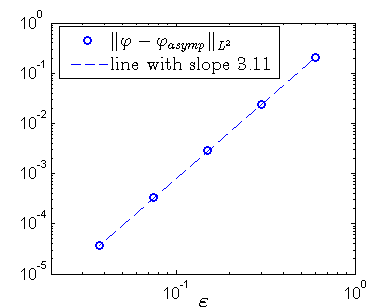}}
  \end{center}
  \caption{{\small Left: Bifurcation diagram in the $(\omega,
    \|\phi\|_{L^2(\P^2)})$-plane for example A: $N=1$,
    $k_*=(\tfrac{1}{2},\tfrac{1}{2})$. Dashed lines: approximation
    $\|\phi\|_{L^2(\P^2)} \sim |A|\sqrt{(\omega-\omega_*)/\Omega}$
    with $\Omega=\sigma=\pm 1$. Curves bifurcating to the left/right
    of $\omega_*$ are for $\sigma =\mp 1$, respectively. The spectrum
    $\sigma(-\Delta +V)$ is plotted on the horizontal axis.  
Right: Error for $\sigma =-1$, where 
    $\phi_{asymp} := \eps A \xi_{1}(x, k_*)$.}
}
  \label{F:bif_diag_N1_kM}
\end{figure}

\begin{figure}[!ht]
  \begin{center}
	\hspace{-0.5cm}\scalebox{0.5}{\includegraphics{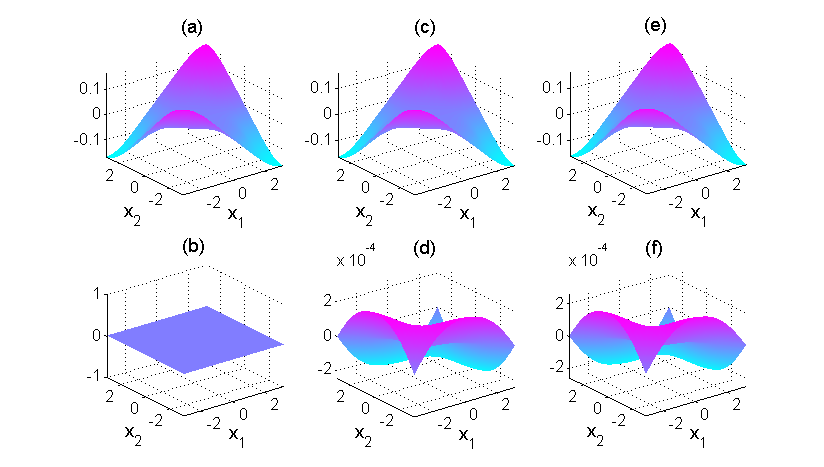}}
\hspace{-1cm}\scalebox{0.5}{\includegraphics{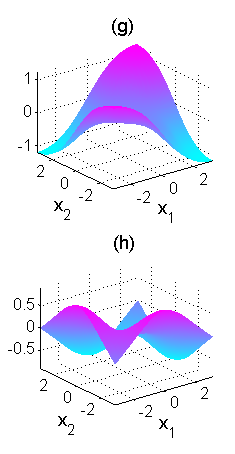}}
  \end{center}
 \caption{{\small Nonlinear Bloch waves for example A. (a) and (b): real and imaginary part of the approximation $\eps A \xi_1(x,(1/2,1/2))$ at $\eps =0.12$;
     (c) and (d): real and imag. part of $\phi$ at (A1-) in
     Fig.~\ref{F:bif_diag_N1_kM}; (e) and (f): real and imag. part of
     $\phi$ at (A1+) ($\sigma
     =1$ and $\omega =\omega_*+\sigma \eps^2$); (g) and (h): 
 real and
     imag. part of $\phi$ at (A2+)}}
  \label{F:prof_N1_kM_near}
\end{figure}

\textbf{Example A.} 
We choose $k_*=(1/2, 1/2)$ and $\pa \sigma(-\Delta +V) \ni
\omega_* = \omega_1(k_*)\approx 1.703$, see point (A) in
Fig.~\ref{F:band_str}.  This leads to $\mu = \|\xi_{1}(\cdot,
k_*)\|_{L^4(\P)}^4\approx 0.0765 $ and choosing $\Omega =\sigma$ and
$\text{arg}(A)=0$, we get
$$A = \tfrac{1}{\sqrt{\mu}}\approx 3.6154.$$
Figure \ref{F:bif_diag_N1_kM} shows the continuation diagram (in the
$(\omega, \|\phi\|_{L^2(\P^2)})$-plane) of the nonlinear Bloch waves
bifurcating from $\omega_*$ for $\sigma=-1$ and $\sigma =1$, 
the asymptotic curves $(\omega_*+\Omega\eps^2, \eps |A|)$ for
$\eps \geq 0$, and the error between the two in the log-log scale. 
The observed convergence rate is 3.11, in 
agreement with Theorem \ref{T:main}.
In Fig.~\ref{F:prof_N1_kM_near} we plot profiles $\phi$ and the
asymptotic approximation $\eps A \xi_{1}(x, k_*)$ at $\omega =\omega_*
+ \eps^2 \Omega$ with $\eps \approx 0.12$, i.e.~close to the
bifurcation point, see points (A1-) and (A1+) in
Fig.~\ref{F:bif_diag_N1_kM}, and $\phi$ at $\omega \approx 2.75$ for
$\sigma =1$, i.e.~far from the bifurcation, cf. point (A2+). 
The asymptotic approximation is real
since the Bloch wave $\xi_{1}(x, k_*)$ has been selected real. This is
possible as $k_*$ is one of the high symmetry points $\Gamma,X,M$. 

\subsection{Numerical Examples for $N=2$.}\label{num2s}

We present computations for two consistent examples with $N=M=2$, cf. \S\ref{S:N2} (a), where the ACMEs \eqref{E:N2_consist} are valid. 
In example B we choose $\omega_* \in \pa
\sigma(-\Delta +V)$ and in example C we take $\omega_* \in
\text{int}(\sigma(-\Delta +V))$.

\textbf{Example B:} 
We choose $k_*^{(1)}=(1/2,0), k_*^{(2)}=(0,1/2), \pa
\sigma(-\Delta +V) \ni \omega_*
=\omega_2(k_*^{(1)})=\omega_2(k_*^{(2)}) \approx 2.035$, see point 
$(B)$ in Fig.~\ref{F:band_str}. Choosing real Bloch waves $\xi_{2}(\cdot, k^{(1)}_*), \xi_{2}(\cdot, k^{(2)}_*)$ (possible due to the real boundary conditions in \eqref{E:Bloch_ev_prob}), we obtain
\begin{align*}
  &\mu_{1111} = \mu_{2222} = \|\xi_{2}(\cdot, k^{(1)}_*)\|_{L^4(\P)}^4\approx  0.0901, \\
  &\mu_{2121}=\mu_{1221}=
  \int_{\P^2}\xi_2(x,k^{(1)}_*)^2\xi_2(x,k^{(2)}_*)^2 dx \approx
  0.003,
\end{align*}
where the equalities between the $\mu$ coefficients follow by the
symmetry $\xi_2((x_1,x_2),k^{(1)}_*)=\xi_2((x_2,x_1),k^{(1)}_*)$ and
the fact that real Bloch waves $\xi_2(x,k^{(1)}_*),
\xi_2(x,k^{(2)}_*)$ have been chosen.

The resulting values of $|A_1|$ and $|A_2|$ are $|A_1|=|A_2|\approx
3.17567$ and in order to satisfy reversibility, we choose zero 
phases, such that
$$A_1=A_2\approx 3.17567.$$
The non-degeneracy condition is satisfied as our computation of the
eigenvalues of ${\bf \hat{J}}$ produces
$$ \lambda_1 \approx -0.1223, \ \lambda_2 =0, \ \lambda_3 \approx 1.6332, \ \lambda_4 =2.$$

\begin{figure}[!ht]
  \begin{center}
    \scalebox{0.5}{\includegraphics{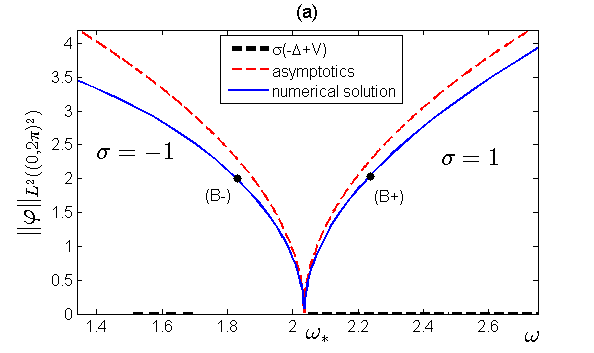}}
  \scalebox{0.53}{\includegraphics{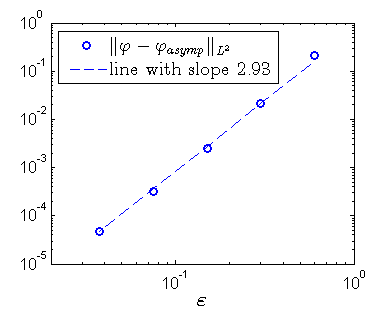}}
  \end{center}
  \caption{{\small Left: Bifurcation diagram in the $(\omega,
    \|\phi\|_{L^2(\P^2)})$-plane for example B: $N=2,
    k_*^{(1)}=(1/2,0), k_*^{(2)}=(0,1/2)$. Full lines: numerically
    computed solution $\phi$; dashed lines: asymptotic approximation
    $\|\phi\|_{L^2(\P^2)} \sim \sqrt{(\omega-\omega_*)/\Omega}
    \|\sum_{j=1}^2A_j\xi_2(\cdot,k^{(j)}_*)\|_{L^2(\P^2)}$ with
    $\Omega=\sigma=\pm 1$. Right:  error for $\sigma =-1$, 
 $\phi_{asymp} :=
    \eps \sum_{j=1}^2A_j \xi_{2}(x, k_*^{(j)})$.}}
  \label{F:bif_diag_N2_kXXprime}
\end{figure}

\begin{figure}[!ht]
  \begin{center}
    \scalebox{0.53}{\includegraphics{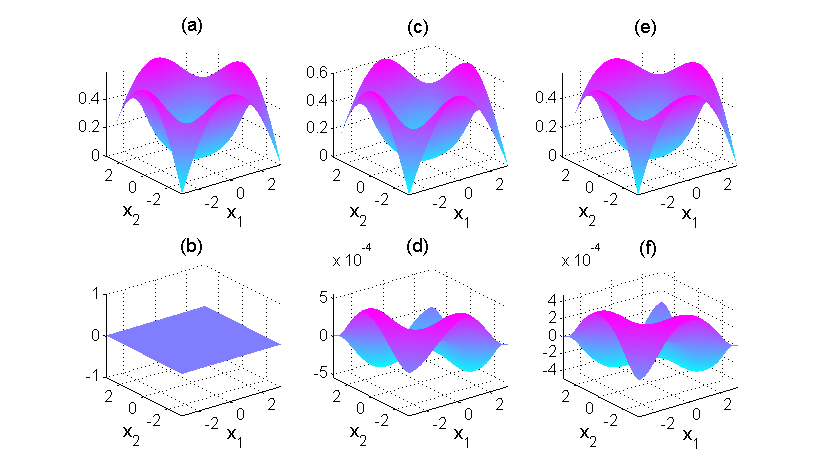}}
  \end{center}
  \caption{{\small Nonlinear Bloch waves for example B. (a) and (b): real and imaginary part
    of the asymptotic approximation $\eps \sum_{j=1}^2A_j \xi_{2}(x,
    k_*^{(j)})$ at $\eps=0.452$; (c) and (d): real and imaginary part
    of $\phi$ at (B-) in Fig.~\ref{F:bif_diag_N2_kXXprime}, i.e. for
    $\sigma =-1$ and $\omega =\omega_*+\sigma \eps^2$; (e) and (f):
    real and imaginary part of $\phi$ at (B+) in
    Fig.~\ref{F:bif_diag_N2_kXXprime}, i.e. for $\sigma =1$ and
    $\omega =\omega_*+\sigma \eps^2$.}}
  \label{F:prof_N2_kXXprime}
\end{figure}

The continuation diagram in Fig.~\ref{F:bif_diag_N2_kXXprime} plots
the families of nonlinear Bloch waves bifurcating from $\omega_*$ for
$\sigma=-1$ and $\sigma =1$, the asymptotic curves
$(\omega_*+\Omega\eps^2, \eps
\|\sum_{j=1}^2A_j\xi_2(\cdot,k^{(j)}_*)\|_{L^2(\P^2)})$ for $\eps \geq
0$, and the $\eps-$convergence
of the approximation error for this case. 
The solutions $\phi$ at the points (B-), i.e. $\omega=1.8304$, and
(B+), i.e. $\omega = 2.2392$, marked in
Fig.~\ref{F:bif_diag_N2_kXXprime} are plotted in
Fig.~\ref{F:prof_N2_kXXprime} together with the asymptotic
approximation $\eps \sum_{j=1}^2A_j \xi_{2}(x, k_*^{(j)})$ at $\omega
=\omega_* + \eps^2 \Omega$ with $\eps \approx 0.452 \approx
\sqrt{\omega_*-1.8304}\approx \sqrt{2.2392-\omega_*}$. Despite the 
large value of $\eps$ the asymptotic approximation is
relatively good.

\textbf{Example C:} 
Finally, we take $k_*^{(1)}=(1/4,1/4), k_*^{(2)}=(-1/4,-1/4),
\text{int}(\sigma(-\Delta +V)) \ni \omega_*
=\omega_1(k_*^{(1)})=\omega_1(k_*^{(2)}) \approx 1.576$, see Point
(C) in Fig.~\ref{F:band_str}. Fixing the free complex phase of the Bloch
waves by setting
$\text{Im}(\xi_1((0,0),k_*^{(1)})=\text{Im}(\xi_1((0,0),k_*^{(2)})=0$,
we obtain
\begin{align*}
  &\mu_{1111} = \mu_{2222} =\mu_{1221}= \|\xi_{1}(\cdot, k^{(1)}_*)\|_{L^4(\P)}^4\approx  0.0526, \\
  &\mu_{2121} =
  \int_{\P^2}\xi_1(x,k_*^{(2)})^2\overline{\xi_1(x,k_*^{(1)})}^2
  dx\approx 0.0412.
\end{align*}
The identities $\mu_{1111} = \mu_{2222} =\mu_{1221}$ follow 
from $\xi_1(x,k_*^{(2)}) = \overline{\xi_1(x,k_*^{(1)})}$, and 
$\mu_{2121} \in \R$ follows because
$\text{Im}(\xi_1(x,k_*^{(1,2)}))$ happen to be antisymmetric in the
$x_1=x_2$ direction. 
The resulting values of $A_1$ and $A_2$ (once again selected real due
to $\mu_{2121} \in \R$) are
$$A_1=A_2\approx 2.242.$$
Also here the non-degeneracy condition is satisfied as our computation
of the eigenvalues of ${\bf \hat{J}}$ produces $\lambda_1 \approx
-0.9427, \ \lambda_2 \approx -0.828, \ \lambda_3 = 0, \ \lambda_4 =
2.$

\begin{figure}[!ht]
  \begin{center}
    \scalebox{0.5}{\includegraphics{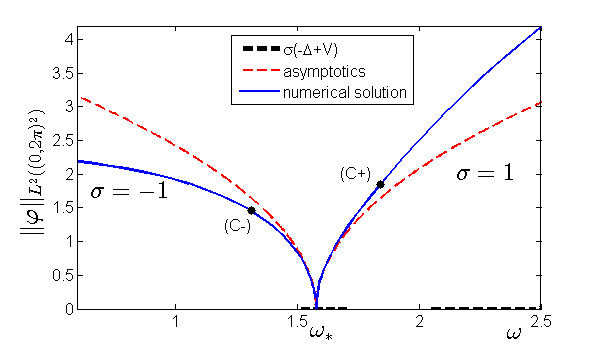}}
 \scalebox{0.53}{\includegraphics{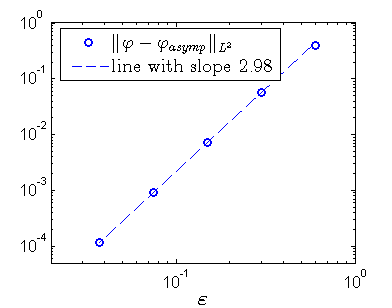}}
  \end{center}
  \caption{{\small Bifurcation diagram in the $(\omega,
    \|\phi\|_{L^2(\P^2)})$-plane and  error for 
$\sigma =-1$ for example C: $N=2,
    k_*^{(1)}=(1/4,1/4), k_*^{(2)}=-k_*^{(1)}$.}}
  \label{F:bif_diag_N2_kdiag}
\end{figure}

\begin{figure}[!ht]
  \begin{center}
    \scalebox{0.53}{\includegraphics{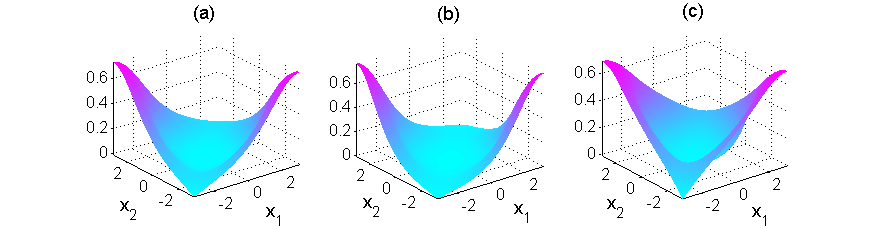}}
  \end{center}
  \caption{{\small Nonlinear Bloch waves  for example C. (a): the approximation $\eps
    \sum_{j=1}^2A_j \xi_{1}(x, k_*^{(j)})$ with $\eps =0.516$; (b) and
    (c): $\phi$ at (C-) and (C+) resp. in
    Fig.~\ref{F:bif_diag_N2_kdiag}. In (b) $\sigma =-1$ and $\omega
    =\omega_*- \eps^2 \approx 1.31$ and in (c) $\sigma =1$ and $\omega
    =\omega_*+ \eps^2 \approx 1.842$.}}
  \label{F:prof_N2_kdiag}
\end{figure}

The continuation diagram from $\omega_*$ for $\sigma=-1$ and $\sigma
=1$ and an error plot for $\sig=-1$ are in Fig.~\ref{F:bif_diag_N2_kdiag}, 
and the solutions $\phi$ at the
points (C$\mp$) with $\omega=1.31, \omega = 1.842$, are in
Fig.~\ref{F:prof_N2_kdiag} together with the asymptotic approximation
$\eps \sum_{j=1}^2A_j \xi_{1}(x, k_*^{(j)})$ at $\omega =\omega_* +
\eps^2 \Omega$ with $\eps \approx 0.516 \approx
\sqrt{\omega_*-1.31}\approx \sqrt{1.842-\omega_*}$.

\section{Gap solitons, out--of--gap solitons, and tNLBs}
\label{ogssec}

NLBs play an important role in the bifurcation structure of many other solutions of \reff{E:SNLS_dD}. As the numerical computations below suggest, when solutions with decaying tails are continued from spectral gaps into spectrum of $-\Delta +V$, they delocalize as the tails become oscillatory with the oscillation structure agreeing with a certain NLB. This puts NLBs in a strong connection with other prominent solutions of \reff{E:SNLS_dD}.

\subsection{1D simulations}\label{1dnum}
\begin{figure}[ht]
  \begin{center}
\ig[width=28mm,height=38mm]{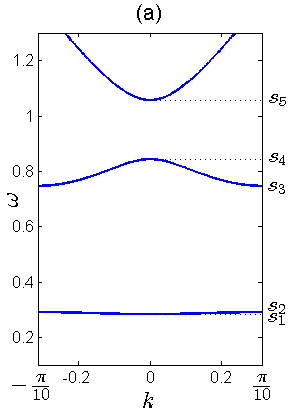}
\ig[width=76mm,height=35mm]{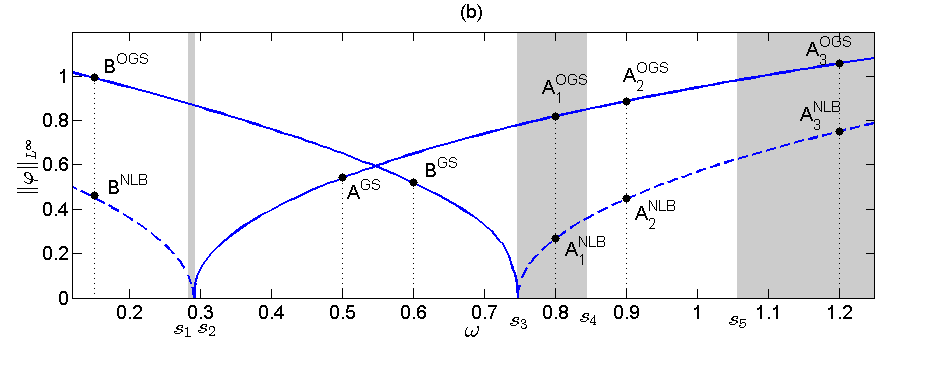}\hspace{1mm}
\raisebox{15mm}{\begin{minipage}{58mm}
\ig[width=29mm,height=30mm]{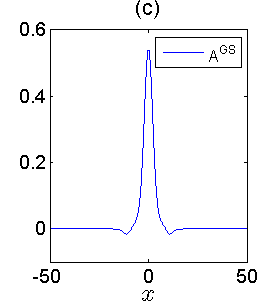}\ig[width=29mm,height=30mm]{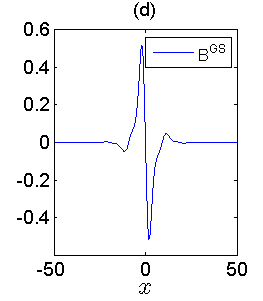}
\end{minipage}}

\ig[width=90mm]{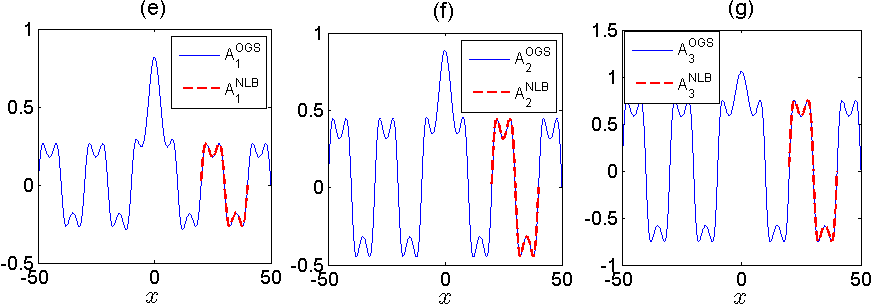}
\ig[width=32mm]{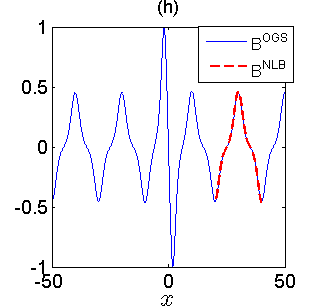}

\end{center}
 \caption{Panel (a): 
the  band--structure for  $V(x)=\sin^2(\frac{\pi x}{10})$; (b): 
a bifurcation diagram of GS and NLBs for \reff{E:SNLS_dD} in 1D, 
$\sigma=1$. 
(c) and (d): plots of the GS at $A^{\text{GS}}$ and $B^{\text{GS}}$ resp.; (e)-(h): OGS and NLB at the remaining marked points.}\label{ogs1Df1}
\end{figure}

We  first consider \reff{E:SNLS_dD} in 1D with $V(x)=\sin^2(\frac{\pi x}{10})$, 
which is a standard choice in 1D. See Fig.~\ref{ogs1Df1}(a)  
for the band--structure, which shows the gaps $(s_2,s_3)$ and 
$(s_4,s_5)$. The first five spectral edges are, approximately, 
$$
s_1\approx 0.2832,\quad s_2\approx 0.2905,\quad s_3\approx 0.7468, \quad 
s_4\approx 0.8434, \quad s_5\approx 1.0568. 
$$
For suitable $\sigma=\pm 1$, 
so called gap solitons bifurcate from the edges into a gap 
\cite{Ac00, Ag01b,dmitryb}.
We display here gap soliton families bifurcating
for $\sigma=1$ to the right from edge $s_2$ and for $\sigma=-1$ to the left from $s_3$. To study these numerically, 
we consider \reff{E:SNLS_dD} on a large domain $x\in(-100,100)$ 
with Neumann boundary conditions, 
and obtain the bifurcation diagram in Fig.~\ref{ogs1Df1}(b), where moreover we restrict to real solutions. 

The gap solitons can be continued in $\omega$ well into the gap. 
In fact, numerically they can also be continued into the next spectral band (and 
even further into higher gaps and bands), where they are 
called out--of-gap solitons (OGS) \cite{ys03, JKKK11}. During this continuation 
the tails of the OGS pick up the oscillations from the NLB that 
bifurcates at the first gap edge, where the continuation family enters the spectrum, i.e., the $A$ branch of OGS 
picks up the oscillations from the NLB branch $A^{\text{NLB}}$ that 
bifurcates at $s_3$ to the right. Moreover, the numerics then 
show that the tails of the OGS are given by $A^{\text{NLB}}$ for \textit{all} $\om>s_3$. The same can be observed for the $B$-family, where the OGS tails
are given by $B^{\text{NLB}}$ for all $\om<s_2$. Thus, an OGS is a homoclinic orbit to a NLB. 

Besides GSs the NLB play a role in the delocalization of many other solutions. In Fig.~\ref{ogs1Df2} we show two other solution branches for illustration. The $B$ branch is an example of a so called {\em truncated NLB} (tNLB), \cite{aok06, 
wyak09, zlw09}. Point $B_0$ at $\om=0.5$ on that branch 
is obtained from using 
\huga{\label{igform} 
\phi_{IG}(x)=a{\rm sech}(x^2/w), \quad a=0.5, w=50, 
}
as an initial guess for a Newton--loop for \reff{E:SNLS_dD}. 
It is homoclinic to $0$ and composed of three periods of the NLB bifurcating from 
$s_1$ in the middle. 
That is why such solutions are called truncated NLBs. 
By varying, e.g., $w$ in \reff{igform}, we can in fact produce 
tNLBs composed of any number of periods of the NLB. 

\begin{figure}[ht]
  \begin{center}
\ig[width=70mm,height=60mm]{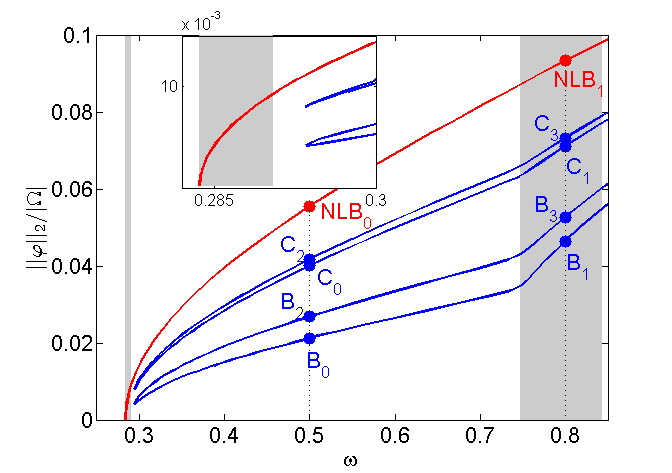}\hspace{7mm}
\raisebox{30mm}{\begin{minipage}{75mm}
\ig[width=35mm,height=30mm]{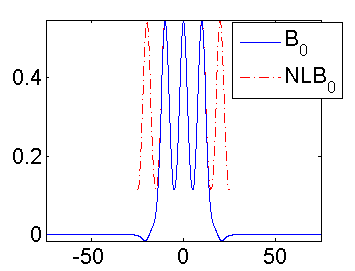}\ig[width=35mm,height=30mm]{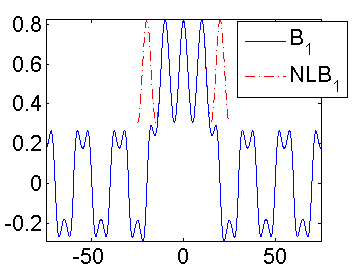}
\ig[width=35mm,height=30mm]{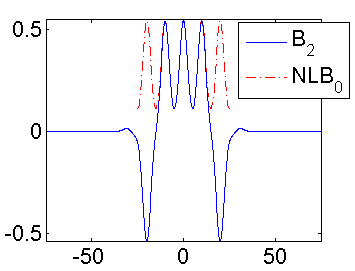}\ig[width=35mm,height=30mm]{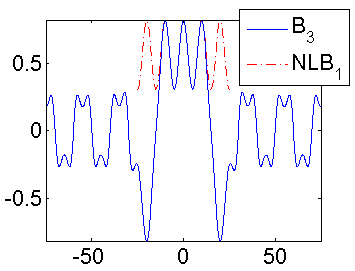}
\end{minipage}}

\ig[width=35mm,height=30mm]{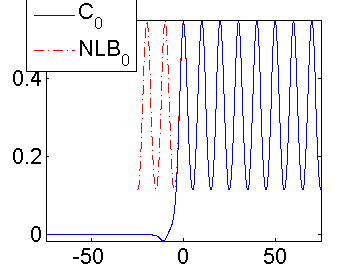}\ig[width=35mm,height=30mm]{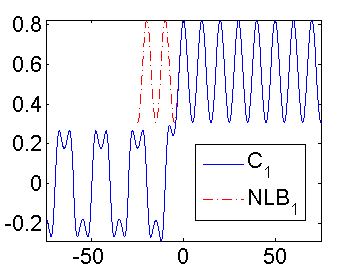}
\ig[width=35mm,height=30mm]{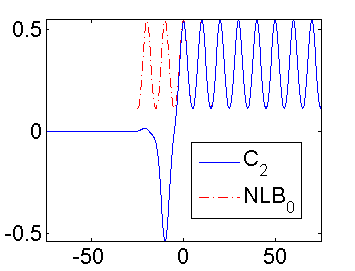}\ig[width=35mm,height=30mm]{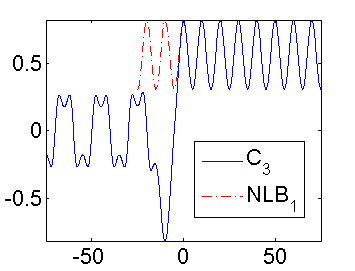}

\end{center}
 \caption{{\small Big panel: bifurcation diagram of first NLB and example 
tNLBs and heteroclinics for  $V(x)=\sin^2(\frac{\pi x}{10})$, $\sig=1$; 
spectral bands in grey. For the $B$ (tNLB) and $C$ (heteroclinic) 
branches we continue from $\om=0.5$ in the positive and negative $\om$ direction.  
For the negative $\om$--directions we obtain folds close to $s_2$, 
cf.~the inset. The NLB branch bifurcates from $s_1$. Smaller panels: example 
plots, where the red dash-dotted line indicates the $s_1$--NLB at the respective 
$\om$--values. The tails of the tNLBs, and the zero--level of the 
heteroclinic, pick up the $s_3$--NLB when entering the second band.}}\label{ogs1Df2}
\end{figure}

An important feature of tNLBs is that they do {\em not} 
bifurcate from $0$, in contrast to the GS. In fact, as a 
tNLB approaches the gap edge next to the $\om$ value where 
its building--block NLB bifurcates, it turns around while 
picking up a negative copy of the pertinent NLB. See also 
\cite{wy08} for a further discussion (in 2D). 
On the other hand, tNLBs behave quite similarly to GS 
upon continuation through the other gap--edge: the tails 
again pick up the NLB bifurcating at the edge (the tNLBs in Fig.~\ref{ogs1Df2} pick up the NLB family bifurcating from $s_3$ in Fig.~\ref{ogs1Df1}), and afterward can be continued to arbitrarily large $\om$ as homoclinics to 
these NLBs, still being close to the original NLB in the middle. For these delocalized tNLBs we suggest the acronym dtNLBs. 

Finally, as there are ``arbitrarily long'' tNLBs, it is not surprising 
that there also exist heteroclinics between $0$ and NLBs. Upon 
continuation in $\om$ these 
essentially behave like tNLBs, see the $C$ branch in Fig.~\ref{ogs1Df2} 
for an example. 

A rigorous analysis of OGS, tNLBs, dtNLBs, 
and the above heteroclinics remains an intriguing open problem, 
even in 1D.  For the 1D case with {\em narrow} gaps 
a system of first order differential coupled mode equations for the 
envelopes of the linear gap edge Bloch waves is derived in \cite{ys03}.
Under suitable conditions, this system has spatial homoclinic orbits  
to nonzero fixed points, which thus corresponds to dtNLBs or OGS. 
However, presently it is not clear how to make this analysis rigorous. 
In \cite{JKKK11} some explicit OGS solutions are given for the case of 
a 1D discrete NLS.  
Concerning tNLBs, \cite{zw09} gives so called composition relations, 
which however are rather heuristic. 
Delocalized (or generalized) solitary waves also occur in other 
nonlinear equations, in particular from fluid dynamics. 
See, e.g., \cite[Chapters 1 and 6.4]{boyd98} for a review, and a guide 
to the literature for rigorous existence proofs, for instance 
\cite{SC94} for the case of the fifth-order KdV equation. The solutions studied in this literature, however, typically have exponentially small tails, which is different from our OGS and dtNLBs, where the amplitude of the tails is that of the NLBs, i.e. $O(\eps^{1/2})$ in the bifurcation parameter $\eps$ as given by \eqref{E:ansatz_mult}.

\subsection{2D simulations}
In 2D similar effects as in Figs.~\ref{ogs1Df1} and \ref{ogs1Df2} 
occur, but the solution structure becomes much richer, also 
qualitatively. For instance, since in 1D the pertinent NLS amplitude equation 
is scalar, there typically is only  one 
 GS bifurcating at some $s_j$ (modulo phase invariance, and 
on--site or off-site effects, see \cite{dmitryb}).  In 2D, in many cases 
the GS are described by {\em systems} of NLS equations, see \cite{DU09}, 
and there may be various different GSs.  Moreover, while typically in 1D 
different tNLBs at fixed 
$\om$ only differ in the number of NLB periods, and the number and arrangements 
of ``ups'' and ``downs'', in 2D we can easily produce qualitatively different 
tNLBs. Accordingly, in the references already cited, in particular 
\cite{wy08}, various families of 2D tNLBs have been studied, with focus on 
the fold structure near one of the gap-edges. 

However, the continuation of either tNLBs or GSs into the other 
spectral band seems to be much less studied, but see also \cite{yang10}. 
Here we 
restrict ourselves to just illustrating the continuation of 
two (real) families of GS to OGS. 
We return to the potential \reff{E:V}, and in Fig.~\ref{ogsf1} 
continue the $\sig=\pm 1$ GS from the first gap into the respective 
other spectral band. Numerically we again use a large domain 
$x\in(-40\pi,40\pi)^2$ 
with Neumann  boundary conditions. 
For the GS these boundary conditions hardly matter, but the 
way in which the tails pick up NLBs as the GS enter the spectral 
bands does significantly depend on the boundary conditions, as should be expected. 
For instance, in (b) we find dislocations in the tail patterns 
along the coordinate axes, and in (e) along one of the diagonals. 
Numerically, these dislocations strongly depend on the chosen 
domain size and boundary conditions. Nevertheless, in all cases considered 
the tails of the GS again pick up a 
pertinent NLB in large parts of the domains. Figure \ref{ogsf1} 
just gives two illustrations. 
\begin{figure}[!ht]
  \begin{center}
\ig[width=120mm]{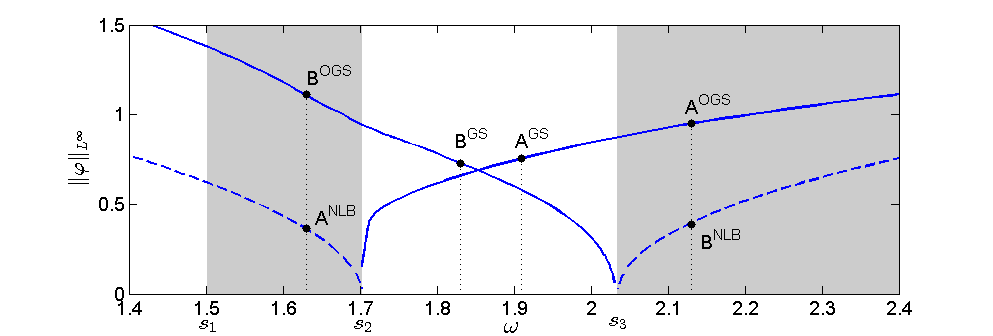}\\
\ig[width=55mm]{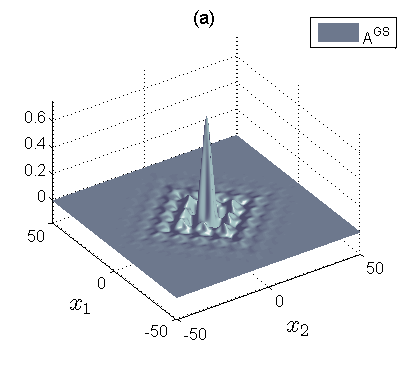}\ig[width=95mm]{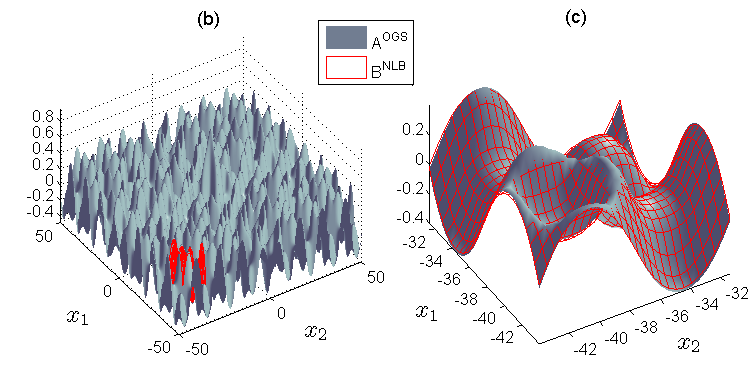}\\
\ig[width=55mm]{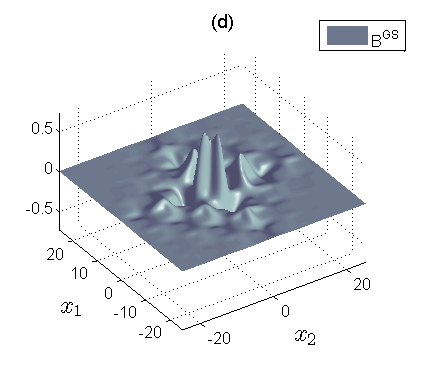}\ig[width=95mm]{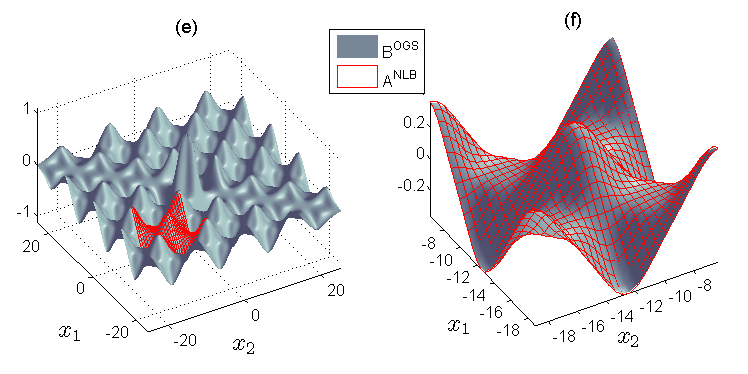}
 \end{center}
 \caption{{\small Continuation of 2D--GS from the first spectral gap 
to OGS, see the text for comments.}\label{ogsf1}}
\end{figure}

Note that in 2D there is typically a number of NLB families
bifurcating from a given point in the spectrum, cf. Theorem
\ref{T:main} with $N>1$. At $\omega=s_3$ the level set of the band
structure is $\{(1/2,0),(0,1/2)\}$, i.e.~to capture at least all the
NLBs predicted by Theorem \ref{T:main} to bifurcate from $s_3$, we
must take $N=2$. The resulting ACMEs are \eqref{E:N2_consist}. The NLB
family $B$ in Fig. \ref{ogsf1}, which happens to describe the
oscillations in $A^{\text{OGS}}$, has $A_1=A_2 \in \R$. In general it
is not clear how to choose the correct solution of ACME which matches
the tail oscillations in a given OGS in $d$D with $d\geq 2$.

\subsection{Remarks on stability of NLBs, GS, tNLBs, and OGS and dtNLBs}
Dynamic stability of GS, NLBs, (localized and delocalized) tNLBs and
OGS is an important but widely open question. Previous work, mostly
based on numerics and formal asymptotics, include the following: In
\cite{HAY2011} and \cite{BD11} it is shown that, so called, on-site
GSs in 1D are spectrally stable while off-site ones are
unstable. \cite{SWCY2008} and \cite{Y2010} show that 2D GSs near the
spectral edge from which they bifurcate, are spectrally unstable but
can be stable further away from the edge.  Next, \cite{GSS12} gives
numerics and formal asymptotics that indicate that GS near the middle
of (narrow) band gaps may be unstable due to a four wave mixing with
the gap edge NLBs.  In \cite{KS02} it is discussed that modulational
instability of 1D NLBs bifurcating from gap edges into the gap can
lead to the formation of GS, while NLBs bifurcating from the edge into
the spectrum are stable.  Some semi-analytical results on the
stability of NLB at the bottom of the band structure are given in
\cite{CP2012}, where, together with the secondary bifurcations from
NLBs, exchange of stability results are derived under some
assumptions, and numerical justifications and comparisons to numerical
time--integration are given.  Regarding tNLBs, in \cite{wyak09} it is
shown numerically that tNLB of the type $B_0$ in Fig.~\ref{ogs1Df2},
i.e., consisting of arbitrary many periods of NLBs of the same parity,
can be stable, while tNLB on the upper branch, consisting of up and
down copies of the basic NLB, are generically unstable.

Here we report some stability results from numerical time 
integration of \reff{tGP} using a Fourier split-step method. We plug 
$\psi(t,x)=\er^{\ri \om t}\phi(t,x)$ into \reff{tGP} to obtain 
\huga{\label{tGP2}
\ri\pa_t\phi=\Delta \phi
-(V(x)-\om) \phi - \sigma |\phi|^2\phi. 
}
Since \eqref{tGP} is Hamiltonian, we can at best obtain spectral 
stability (but not linearized stability) from linearization of 
\reff{tGP2} around a steady state $\phi(x,t)=\phi_0(x)$ of interest (NLB, GS, tNLB etc). 
As initial conditions we choose random perturbations of the steady 
state $\phi_0$, 
with perturbation amplitude 0.1 relative to the amplitude of $\phi_0$. 
These perturbations yield a phase evolution, and thus here we concentrate
 on the solution shape, i.e., we plot the supremum error 
$\||\phi(\cdot,t)|-|\phi_0(\cdot)|\|_{\infty}$ of the modulus. 
This should 
provide an indication regarding stability of the solutions at hand 
and motivate further stability studies.
 Our numerical accuracy was 
checked by using smaller time steps without visible changes of 
results, by comparison with a semi--implicit time stepping, 
and we checked that $\|\phi\|_{L^2}$ was conserved up to 6 digits 
over the rather long time intervals needed in some cases to detect 
instabilities. 

We mostly focus on 1D, and start with the NLBs. In order to draw a connection between stability of NLBs and OGSs, we restrict ourselves to NLBs bifurcating at gap edges in Fig.~\ref{ogs1Df1}.
We choose the periodicity cell 
$x\in\Omega=(-10,10)$ for $\phi$, which corresponds to wave-vectors $k=0$ and $k=1/2$. However, the results 
appear to be the same for larger domains, i.e., if we observe 
instabilities, then they are w.r.t. the same spatial period. 
Below we use the symbol $s_n\pm$ to denote the NLB families 
bifurcating from a spectral edge $s_n$ to the right and left respectively.

\begin{figure}[ht]
  \begin{center}
{\small
\begin{tabular}{cccc}
(a) NLB $s_3+, \om{=}0.76$ &(b) NLB $s_3+, \om{=}0.76$ &(c) NLB $s_3+, \om{=}0.9$ & (d) NLB $s_3+, \om{=}0.9$ \\
$|\phi|$& amplitude error &$|\phi|$& amplitude error \\
\ig[width=40mm,height=36mm]{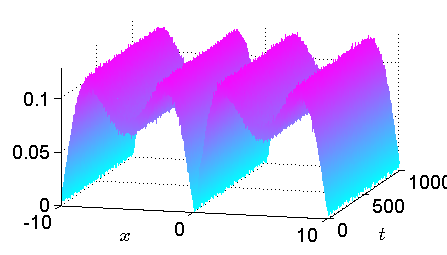}&
\ig[width=35mm,height=30mm]{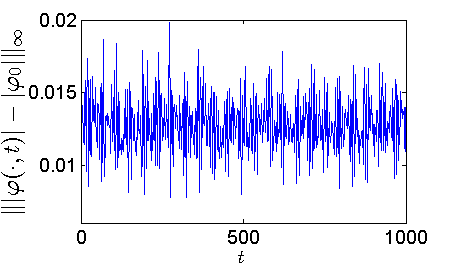}&
\ig[width=40mm,height=36mm]{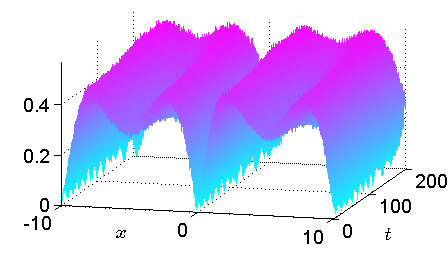}&
\ig[width=35mm,height=30mm]{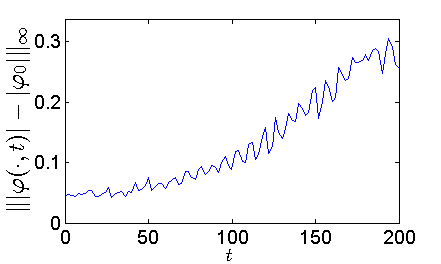}\\
\end{tabular}
\begin{tabular}{cccc} 
(e) OGS $s_2+,\ \om{=}0.76$ & (f) OGS $s_2+,\ \om{=}0.76$ & (g) OGS $s_2+,\ \om{=}0.9$ & (h) OGS $s_2+,\ \om{=}0.9$  \\
$|\phi|$ & amplitude error& $|\phi|$ & amplitude error\\
\ig[width=40mm,height=36mm]{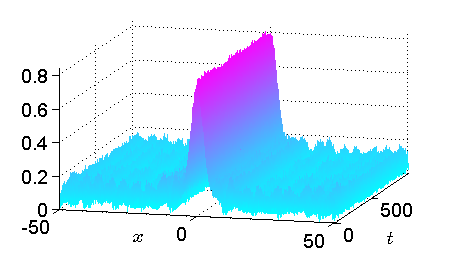}&
\ig[width=35mm,height=30mm]{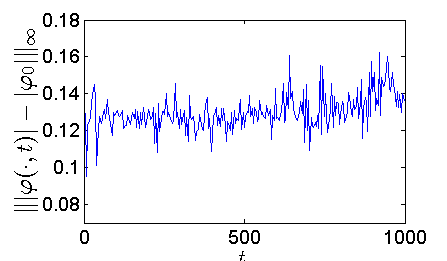}&
\ig[width=40mm,height=36mm]{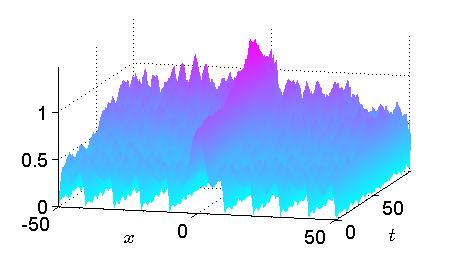}&
\ig[width=35mm,height=30mm]{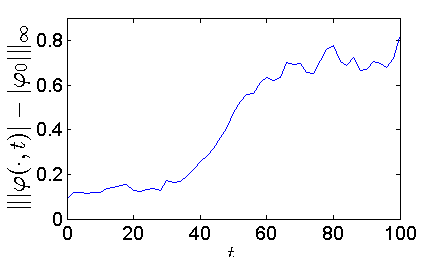}\\
\end{tabular}
}
\end{center}
 \caption{{\small  Numerical integration of \reff{tGP2} with initial data as 
random 10\% perturbations of the indicated NLB in (a)-(d) and OGS in (e)-(h). The NLB-branch $s_3+$ appears stable 
at bifurcation and unstable otherwise. The OGS family $s_2+$ entering the spectrum at $\om =s_3$ seems to inherit the (in)stability of the NLB.}}\label{sf1}
\end{figure}

Panels (a)-(d) in Fig.~\ref{sf1} are for the NLB branch $s_3+$ in
Fig.~\ref{ogs1Df1} for $\sig=1$. As an example in (a),(b) we choose
$\om=0.76$ and observe a stable evolution up to $t=1000$.  In (c), (d)
the NLB at $\om=0.9$, i.e. further away from the bifurcation edge, is
shown to be unstable. An analogous situation occurs for the $s_2-$
family of NLB for $\sigma=-1$. It appears stable for $\om \in
(0.2,s_2)$ and unstable for $\om<0.2$.

Our tests suggest that similar stability results (stability near the
bifurcation edge and instability otherwise) hold also for other NLB
branches bifurcating into the spectral bands 
in Fig.~\ref{ogs1Df1}. The stability near the bifurcation
edge 
is in agreement with
\cite{KS02}. The loss of stability of NLB away from a neighborhood of
the bifurcation point is presumably due to a secondary ``loop''
bifurcation as discussed in \cite{CP2012}. 

Given the results on the NLBs from Fig.~\ref{sf1} (a)-(d), we may expect
that the OGS with tails picking up the $s_3+$ NLBs at $s_3$ 
inherit the (in)stability from the NLBs $s_3+$ NLBs and hence 
is stable at $\om =0.76$ and unstable at $\om=0.9$. 
This is, indeed, observed in Fig.~\ref{sf1}
(e)-(h). Note that one does not expect any instability from the GS
``component'' of the OGS since we work with on-site GSs, which have
been reported in \cite{HAY2011} and \cite{BD11} to be stable. 

Next, the NLB branch $s_1+$ bifurcating from 
$s_{1}$ to the right for $\sig=1$ is stable throughout the first band and 
the first gap. Thus, another relevant question is whether 
the associated tNLBs inherit this stability of their building blocks. 
In accordance with, e.g., \cite{wyak09}, this is the case for tNLB consisting 
of copies of NLBs with the right parity, i.e., only up or only down 
copies of $s_1+$ NLBs; see for instance the tNLBs $B_0$ and 
$C_0$ from Fig.~\ref{ogs1Df2}. In a next step we then studied the stability 
of dtNLBs obtained from the continuation of such tNLBs across the $s_3$ 
gap edge. 
This is in complete agreement with the stability of OGS, i.e.,  
the dtNLBs obtained from $B_0, C_0$ (see $B_1, C_1$ in  Fig.~\ref{ogs1Df2}) 
are stable for small $\om{-}s_3{>}0$ but become unstable 
for larger $\om{-}s_3$. 
On the other hand, we found that tNLBs consisting of up and down copies 
of NLBs (e.g., $B_2$, $C_2$ in Fig.~\ref{ogs1Df2}) are unstable, as in 
\cite{wyak09}. For instance,  for $C_2$ the leading down NLB 
is first converted into an up NLB, and then a defect wanders to the right.

In 2D, the NLBs from Fig.~\ref{f2} all appear modulationally unstable, with 
however very long transients before the instability sets in 
for the branches bifurcating at smaller $\om$, and this also 
holds for the other NLBs, for instance given in \S\ref{num2s}. 
Moreover, the 2D-GS are expected to be unstable near bifurcation, but
may become stable in the middle of gaps,  see
\cite{IW10} for rigorous results in the semi-infinite gap, 
and \cite[\S6.4]{Y2010} for further
heuristics. This agrees with our numerics, where e.g., solutions on 
the $A^{\text{GS}}$ branch from Fig.~\ref{ogsf1} are 
numerically unstable for $\om<1.8$, then 
stable up to $\om=s_3$. However, the OGS for 
$\om>s_3$ with tails containing $B^{\text{NLB}}$ NLBs is clearly 
unstable numerically.    

Thus,  besides analytical results regarding the existence and stability 
of tNLBs, OGS and (d)tNLB, an interesting open problem
 is whether in 2D there exist potentials V such that 
\bcen 
\item \reff{tGP} has stable GS; 
\item \reff{tGP} has stable NLBs bifurcating from gap edges; 
\item putting (1) and (2) together: \reff{tGP} has stable OGS. 
\ecen

\section*{Acknowledgments}
The authors thank Michael I. Weinstein for fruitful discussions, in particular for inquiring about the possibility to generalize the bifurcation assumptions to multiple Bloch eigenvalues, as formulated in (H3). The research of T.D. is partly supported by the \emph{German Research Foundation}, DFG grant No. DO1467/3-1.



\bibliographystyle{plain}
\bibliography{biblio_NLB}

\end{document}